\documentclass[english,letterpaper,12pt,leqno]{article}
\usepackage{hyperref}

\usepackage{tikz}
\usetikzlibrary{decorations.markings}
\usepackage{fullpage}

\usepackage{amsxtra}
\usepackage{amsmath}
\usepackage{amssymb}
\usepackage{amsfonts}
\usepackage{amsthm}

\usepackage[all,2cell]{xy}
\UseAllTwocells
\usepackage{mathrsfs}
\usepackage{enumitem}
\usepackage{subfigure}
\usepackage{float}
 \usepackage{caption}

\newtheorem{thm}[equation]{Theorem}
\newtheorem{cor}[equation]{Corollary}
\newtheorem{lem}[equation]{Lemma}
\newtheorem{prop}[equation]{Proposition}

\newtheoremstyle{example}{\topsep}{\topsep}%
{}%         Body font
{}%         Indent amount (empty = no indent, \parindent = para indent)
{\bfseries}% Thm head font
{.}%        Punctuation after thm head
{2pt}%     Space after thm head (\newline = linebreak)
{\thmname{#1}\thmnumber{ #2}\thmnote{ #3}}%         Thm head spec

\theoremstyle{example}

\newtheorem{Defi}[equation]{Definition}
\newtheorem{defi}[equation]{Definition}
\newtheorem{rem}[equation]{Remark}
\newtheorem{rems}[equation]{Remarks}

\newtheorem{ex}[equation]{Example}

\newtheorem{que}[equation]{Question}

\newtheorem{exa}[equation]{Example}

\newcounter{sparagraphcounter}[subsection]

\def\sparagraph#1{\refstepcounter{sparagraphcounter}
\vskip 1em%
\noindent {\bf \Alph{sparagraphcounter}. #1}%
}

\def\Bc{\mathcal {B}}

\def\CC{\mathbb{C}}
\def\ZZ{\mathbb{Z}}

\def\H{\mathrm{H}}
\def\C{\mathrm{C}}

\def\C{\mathcal{C}}

\def\Cc{\mathca{C}}

\def\p{\mathfrak{p}}
\def\sl{\mathfrak{sl}}
\def\gl{\mathfrak{gl}}

\def\g{\mathfrak{g}}
\def\m{\mathfrak{m}}
\def\b{\mathfrak{b}}
\def\h{\mathfrak{h}}

\def\u{\mathfrak{u}}
\def\p{\mathfrak{p}}

\def\L{\mathfrak{L}}
\def\Sen{\mathfrak{S}}

\def\n{\mathfrak{n}}

\def\ken {\mathfrak{k}}

\def\ZZ{\mathbb{Z}}

\def\RR{\mathbb{R}}
\def\CC{\mathbb{C}}
\def\GG{\mathbb{G}}

\def\C{\mathrm{C}}

\def\Mc{\mathcal{M}}
\def\Oc{\mathcal{O}}

\def\be{\begin{equation}}
\def\ee{\end{equation}}

\def\ev{{\bar 0}}
\def\od{{\bar 1}}
\def\ch{\on{char}}
\def\Com{\on{Com}}
\def\Comc{\mathcal{C}om}
\def\Coker{\on{Coker}}

\def\deg{\on{deg}}
\def\dgAlg{\on{dgAlg}}

\def\gen{\mathfrak{g}}

\def\Hom{\on{Hom}}
\def\Ho{\on{Ho}}
\def\HC{\on{HC}}
\def\det{\text{\rm det}}

\def\Ext{\mathop{\mathrm{Ext}}\nolimits}

\def\from{\colon}
\def\I{\mathrm{I}}
\def\II{\mathrm{I\!I}}
\def\III{\mathrm{I\!I\!I}}

\def\isom{\simeq}

\def\Ker{\mathop{\mathrm{Ker}}\nolimits}
\def\k{\mathbf{k}}

\def\Lie{{\on{Lie}}}
\def\lra{\longrightarrow}

\def\mls#1{\mathop{<}\limits_{#1}}
\def\modd{{m_{\bar 1}}}
\def\mev{{m_{\bar 0}}}

\def\ol{\overline}
\def\on{\operatorname}

\def\pt{\on{pt}}

\def\rev{r_{\ev}}
\def\res{\rho}
\def\rod{r_{{\od}}}

\def\Spec{\on{Spec}}

\def\Vect{\on{Vect}}

\def\wt\widetilde

\def\A{\mathcal A}

\def\Cc{\mathcal {C}}

\def\on{\operatorname}

\newlist{exaenumerate}{enumerate}{1}
\setlist[exaenumerate]{topsep=0pc,itemindent=2.5pc,leftmargin=0pc,label=(\alph{*}),parsep=0pc,itemsep=0pc}

\numberwithin{paragraph}{subsection}

\usepackage{color}

  \newenvironment{dedication}
        {\vspace{3ex}\begin{quotation}\begin{center}\begin{em}}
        {\par\end{em}\end{center}\end{quotation}\vspace{5ex}}

\begin{document} 

\title{Derived varieties of complexes and Kostant's theorem for $\gl(m|n)$. }
\author{ M. Kapranov, S. Pimenov}
\maketitle

\begin{dedication} To Maxim Kontsevich for his 50th birthday
\end{dedication}

\tableofcontents

\vfill\eject

\numberwithin{equation}{section}

\addtocounter{section}{-1}

\section{Introduction} 

\sparagraph{The object of study.}
 Let $\k$ be a field and
 $V^\bullet$ be a finite-dimensional graded 
$\k$-vector space.  The {\em Buchsbaum-Eisenbud variety of complexes} \cite{kempf} associated to $V^\bullet$
is the scheme formed by 
differentials $D$ making $(V^\bullet, D)$ into a cochain complex.
In other words, it
is the closed subscheme
\[
\on{Com}(V^\bullet) \,=\, \biggl\{ D=(D_i) \in \prod_i \,\, {\Hom(V^i, V^{i+1})}\,\, \biggl| \,\, D_{i+1}\circ D_i =0,\,
\forall i \biggr\}.
\]

\vskip .2cm

The group $GL(V^\bullet) = \prod_i GL(V^i)$ acts on $\Com(V^\bullet)$ by automorphisms.
The scheme $\Com(V^\bullet)$ has been studied in several papers \cite{kempf} \cite{deconcini-strickland}
\cite{musili-seshadri}
\cite{lakshmibai}
\cite{gonciulea} \cite{cukierman} where, in particular,  the following remarkable properties
have been found:

\begin{enumerate}
\item[(1)]  $\Com(V^\bullet)$ is reduced (i.e., is an affine algebraic variety), 
and each its irreducible component is normal. 

\item[(2)] If $\ch(\k)=0$, then the coordinate algebra $\k[\Com(V^\bullet)]$ has, as a $GL(V^\bullet)$-module, simple spectrum:
each irreducible representation appears at most once, and one has an explicit description of the
sequences $[\alpha] = (\alpha^{(1)}, ..., \alpha^{(n)})$ of dominant weights for the $GL(V^i)$ for which
the corresponding representation does appear.  More generally, for any $\k$ the algebra 
$\k[\Com(V^\bullet)]$ has a distinguished basis and a straightening rule. 
\end{enumerate}
By definition, $\Com(V^\bullet)$ is cut out, in the affine space 
$ \prod_i \, {\Hom(V^i, V^{i+1})}$,  by a system of homogeneous quadratic equations 
corresponding to the matrix elements of the compositions $D_{i+1}\circ D_i$. 
The complexity of $\Com(V^\bullet)$ stems from the fact that these equations are far from being
``independent": we do not have a complete intersection. 
The goal of the present paper is to study a natural {\em derived extension} of $\Com(V^\bullet)$
which is a {\em derived scheme}  $R\Com(V^\bullet)$, the spectrum of a certain commutative dg-algebra
$A^\bullet(V^\bullet)$. The dg-algebra $A^\bullet(V^\bullet)$
is obtained 
by complementing the above quadratic equations by natural relations among them, then by natural syzygies 
and so on, see \S \ref{subsec:def-init}. Alternatively, we can look at $R\Com(V^\bullet)$,
in terms of the functor it represents (on the category of dg-algebras). From this point of view,  
points of $R\Com(V^\bullet)$ correspond to {\em twisted complexes} 
\cite{BK}\cite{FSbook}, see Proposition \ref {prop:twis-com}. 
 We show that some of the classical properties of $\Com(V^\bullet)$
are extended to $R\Com(V^\bullet)$. 

There are several  additional reasons for interest in $R\Com(V^\bullet)$.

\sparagraph{Derived moduli spaces.} First, $R\Com(V^\bullet)$ can be seen as a typical example of a {\em derived moduli space.}
The program of derived deformation theory initiated by V. Drinfeld \cite{Dr} and M. Kontsevich \cite{kontsevich-course} \cite{kontsevich-torus},
proposes that any  moduli space (i.e., moduli scheme or moduli stack) $\Mc$ appearing in algebraic geometry 
is really a shadow of
a more fundamental derived object $R\Mc$ which is a derived scheme or derived stack,
in particular, (local) functions on $R\Mc$ form a sheaf of  commutative  dg-algebras $(\Oc^\bullet_{R\Mc},d)$
(in the case $\ch(\k)>0$ one should work with simplicial algebras). A consequence of this point of
view for classical algebraic geometry is that $\Mc$ 
(any classical moduli space!) carries  a natural  ``higher structure sheaf"
which is the sheaf of graded commutative algebras $\Oc_\Mc^\bullet = H^\bullet_d(\Oc^\bullet_{R\Mc})$
whose degree 0 component is $\Oc_\Mc$. 

The first step of this program, i.e., a construction of the $R\Mc$ in the first place, is by now largely
accomplished, both in terms of defining the dg-algebras of functions \cite{CFK1}  as well as describing the functors,
represented by the $R\Mc$ on appropriate (model or $\infty$-) categories \cite{lurie-def} \cite{TVe} \cite{Pri},
see \cite{Pri2} for comparison of different approaches.

 However, the geometric significance
(indeed  the very nature) of the sheaves $\Oc^\bullet_\Mc$ remains mysterious. 
Only in the simplest case of {\em quasi-smooth} derived schemes (called $[0,1]$-manifolds in \cite{CFK2})
we have a meaningful interpretation, suggested already in \cite{kontsevich-course}. In this case $\Oc^\bullet_\Mc$ is situated
in finitely many degrees (as $\Oc^\bullet_{R\Mc}$ is essentially a Koszul complex), and the  class
\[
[\Mc]_{\on{virt}}^K \,\,=\,\,\sum_i (-1)^i \on{cl}(\Oc^i_\Mc) \,\,\in \,\, K_0(\Mc)
\]
is the K-theoretic analog of the virtual fundamental class of $\Mc$, see \cite{CFK2} for a precise formulation and proof.

Many of the equations defining  classical moduli spaces can be reduced to  ``zero curvature conditions",
i.e., to some forms of the equation $D^2=0$ defining $\Com(V^\bullet)$ (see \cite{cukierman} for an example
of such reduction).  
The derived schemes $R\Com(V^\bullet)$ provide, therefore,  stepping stones for access to more complicated derived moduli
spaces. In fact, a closely related object, the derived stack $\on{Perf}_\k$ (of perfect complexes
of $\k$-vector spaces), is used by B. Toen and M. Vaqui\'e  as an inductive step for constructing many such spaces in 
 \cite{TVa}. 

From this point of view, our results are suggestive of a general pattern  that should hold for higher structure sheaves 
of many classical moduli spaces. That is, they suggest the following:
\begin{itemize}

\item For moduli spaces $\Mc$ related to higher-dimensional problems (i.e., beyond the quasi-smooth situation),
the sheaf of rings $\Oc^\bullet_\Mc$ is typically not finitely generated over $\Oc_\Mc = \Oc_\Mc^0$.
See simple examples in \S \ref{subsec:def-init}. 

\item Nevertheless, the individual graded components $\Oc_\Mc^i$ are finitely generated (coherent)
over $\Oc_\Mc$, and one should be able to make sense of the statement that the
 generating series 
 \be\label{eq:ser-phi}
F(t) \,\,=\,\,  \sum_i t^i  \on{cl}(\Oc^i_\Mc) \,\,\in \,\, K_0(\Mc)[[t]]
 \ee
 is  ``rational". For   $\Mc=\Com(V^\bullet)$ being the derived variety of complexes, we study, in Chapter \ref {sec:euler},
 a similar generating series with $K_0$ being understood as
 the representation ring of  $GL(V^\bullet)$.

\end{itemize}

\sparagraph {Kostant's theorem. }
 The dg-algebra $A^\bullet(V^\bullet)$ can be defined as the cochain algebra of a certain graded Lie algebra
$\n^\bullet(V^\bullet)$, see \S\ref{sec:def.RCom}. If we contract the $\ZZ$-grading to the $\ZZ/2$-grading,
then the Lie superalgebras $\n^\bullet(V^\bullet)$ for all possible $V^\bullet$ are precisely the nilpotent radicals of
all possible ``Levi-even" parabolic subalgebras in the Lie superalgebras $\mathfrak{gl}(m|n)$
(we refer to \cite{M} for a discussion of Borel and parabolic subalgebras in this setting). 
For a parabolic $\mathfrak p$ in $\mathfrak{gl}(n)$, indeed in any reductive Lie algebra, 
the cohomology of the nilpotent radical $\mathfrak{n}\subset\mathfrak p$ is found by the
classical theorem of Kostant \cite{K}. In particular, this cohomology carries the action of the Levi subgroup $M$
and each representation appears no more than once, in complete parallel with the property (2) above. 

Our main result, Theorem \ref{thm:simple} states that the cohomology of any $\n^\bullet(V^\bullet)$  with constant
coefficients has
a simple spectrum, as a $GL(V^\bullet)$-module. Equivalently the cohomology of the coordinate ring of the derived
variety of complexes has a simple spectrum. 
 This statement can be seen as generalizing both  Kostant's theorem and
the properties  of the varieties of complexes listed above.

The question of generalizing Kostant's theorem to $\mathfrak{gl}(m|n)$ and other reductive Lie
superalgebras has attracted some attention \cite{brundan-stroppel}  \cite{CKL} \cite{CKW} \cite{CZ} \cite{Cou}, 
with the results so far
being restricted to some particular classes of parabolics but treating more general coefficients. 
 Unlike the purely even case, the problem of finding the $\n$-cohomology of arbitrary finite-dimensional
irreducible modules $L_\mu$ of $\gl(m|n)$ is very complicated.
 In fact,  the generating functions for $GL(V^\bullet)$-isotypical components 
of $H^\bullet(\n, L_\mu)$
are Kazhdan-Lusztig polynomials of type $A$
(Brundan's conjecture \cite{Br1} 
 proved in \cite{CLW} with some features interpreted earlier in  \cite{Br2}\cite{CWZ}). 
  Our result for trivial coefficients $\mu =0$ (and arbitrary parabolics)
 means that the correspondng KL  polynomials are monomials with coefficient 1. 
 We do not know if this can be seen from the geometry of Schubert varieties.

 \sparagraph{Structure of the paper and the method of the proof.}
 Chapter \ref{sec:def} defines the derived varieties of complexes, and puts them into two frameworks:
the representation-theoretic framework similar to that of  \cite{kempf} \cite{deconcini-strickland}
and the Lie superalgebra framework, where various degrees are mixed together into one $\ZZ/2$-parity. 

\vskip .2cm

In Chapter \ref{sec:euler} we study the cohomology of $A^\bullet(V^\bullet)$,
 the coordinate dg-algebra of $R\Com(V^\bullet)$
(in the Lie superalgebra framework) at the level of Euler characteristics,
proving that the Euler characteristic of each isotypic component of the cohomology
is either $0$ or $\pm 1$. This is done for several reasons.

First, the evaluation of the Euler characteristic provides a logical intermediary step for our proof of Theorem 
\ref{thm:simple}.

Second,  the proof at the level of  generating functions is easy to understand conceptually.
That is,  in the case of the standard Borel subgroup, the Kostant cohomology with coefficients in $\k$ is easily found.
This is based on the  Cauchy formula  
\be\label{eq:CF}
S^p(V\otimes W) \,\,\simeq \,\, \bigoplus_{|\alpha|=p} \Sigma^\alpha(V) \otimes \Sigma^\alpha(W)
\ee
(simple spectrum decomposition)
and the classical Kostant theorem with coefficients. This leads to a simple formula for the generating
series as a rational function.
The change of the Borel subgroup amounts to re-expansion of this function in different regions which is then
analyzed in terms of adding formal series for  the Dirac $\delta$-function times a function of smaller number of variables. 
In this way we get an explicit formula for the coefficients of the generating series
(Theorem \ref{thm:b-expl}) which, after Theorem \ref{thm:simple} is established, becomes a complete
description of the cohomology of $\n$ (i.e.,  of the cohomology of the dg-algebra $A^\bullet(V^\bullet)$).

Third, the simple (rational)  nature of generating functions here gives a reason to expect good rationality properties of the series
\eqref{eq:ser-phi} in more general algebro-geometric settings. 

\vskip .2cm

One can view our proof of Theorem  \ref{thm:simple} as a ``categorification" of the generating function argument of Chapter 
 \ref{sec:euler}. In particular, the role of the $\delta$-function is played by the Lie superalgebra $\sl(1|1)$
 which appears naturally in the context of an ``odd reflection". In fact, we use some rather detailed information
 about representations and cohomology of $\sl(1|1)$. This information is recalled in Chapter \ref{sec:rec}.
 Another  subject that is recalled is that of {\em mixed complexes}, the algebraic analog of topological spaces
 with circle action, and especially the analog of the concept of {\em equivariant formality} of \cite{EF}. 
 This formality concept allows us, sometimes,  to identify the classes, in the representation ring, of the cohomology
 of two equivariant complexes without these complexes being quasi-isomorphic. 
 
 \vskip .2cm
 
 The proof of Theorem  \ref{thm:simple} is given in Chapter \ref{sec:proof-ss}. It is organized
 in a way parallel to Chapter
 \ref{sec:euler} so that one can see how various simple steps involving generating functions are upgraded to
  spectral sequences, equivariant formality and so on.

\vskip .3cm

\sparagraph{ Self-dual complexes and other classical superalgebras.}
An interesting version of the variety $\Com(V^\bullet)$ is obtained by considering {\em self-dual complexes}
(also known as {\em algebraic Poincar\'e complexes}, see \cite{APC}). For this, we assume that $V^\bullet$ is equipped with
a nondegenerate symmetric or antisymmetric bilinear form  $\beta$ of some degree $n$, so $V^i$  is identified 
with $(V^{n-i})^*$. 
We can then consider the variety  $\Com(V^\bullet,\beta)$ formed by systems of differentials $D=(D_i)$ such that $D_i=D_{n-1-i}^*$. 
This variety (and its natural derived analog) is similarly related to the Kostant cohomology (i.e., cohomology
of nilpotent radical of parabolics) of 
other classical   (symplectic, orthogonal) Lie superalgebras. It is natural to ask whether our statement about simple spectrum
extend to the Kostant cohomology with coefficients in $\k$
for other  simple Lie superalgebras. For the classical case it may be approached by a method similar to ours,   starting with the
well known analogs of \eqref{eq:CF}
which decompose $S^p(S^2(V))$ and $S^p(\Lambda^2(V))$ as $GL(V)$-modules
with simple spectrum.

\vskip .3cm

We would like to thank J. Brundan for looking at the preliminary version of this paper and 
pointing out the relation of our results to the  Kazhdan-Lusztig theory of \cite {Br1, Br2, CLW}, 
  carefully explaining some details. 
  We also thank  C. Stroppel for  bringing some additional references to our attention
  and K. Coulembier for pointing out an incorrect formulation of Proposition \ref{prop:mod-sl11}
  in an earlier version. 
This work was supported by the  World Premier International Research Center Initiative (WPI), MEXT, Japan
and (M.K.)   by the Max-Planck Institute  f\"ur Mathematik, Bonn, Germany.  The first named author
would like to thank the MPI for hospitality, support and excellent working conditions. 

\vfill\eject

\numberwithin{equation}{subsection}

\section{Definitions and framework. }\label{sec:def.RCom}\label{sec:def}

\subsection
{The definition  of $R\Com(V^\bullet)$ and initial examples.}\label{subsec:def-init}

We assume $\ch(\k)=0$ and denote by  $\Comc_\k$ the category of cochain complexes $X^\bullet$ of
$\k$-vector spaces with its  standard symmetric monoidal structure, given by the usual tensor product of complexes
and the Koszul sign rule for the symmetry
\[
X^\bullet\otimes Y^\bullet\to Y^\bullet\otimes X^\bullet, \quad x\otimes y \mapsto (-1)^{\deg(x) \deg(y)} y\otimes x.
\]
By a {\em dg-category} we mean a category $\Cc$ enriched in $(\Comc_\k, \otimes)$, so for each two objects
$a,b\in\Cc$ we have a complex $\Hom^\bullet(a,b)$ with composition maps satisfying the graded Leibniz rule.
In particular, a dg-category with one object $\pt$
is the same as an associative unital dg-algebra $R^\bullet = \Hom^\bullet (\pt,\pt)$
over $\k$. 
We will be particularly interested in (graded) commutative, unital  dg-algebras over $\k$ 
and denote by $\dgAlg_\k$ the category of such dg-algebras concentrated in degrees $\leq 0$.

\vskip .2cm

Fix a graded  $\k$-vector space $V^\bullet$ as in the Introduction. 
Consider the graded  vector space
\[
\n^\bullet \,\,=\,\,\n^\bullet(V^\bullet) \,\,=\,\,\bigoplus_{i<j} \Hom(V^i, V^j)[i-j]
\]
(i.e., $ \Hom(V^i, V^j)$ is put in the degree $j-i>0$). It is a graded associative algebra
(without unit)
and we consider it as a graded Lie algebra via the supercommutator 
\[
[x,y]= xy- (-1)^{\deg(x)\deg(y)}yx.
\]
Let 
\[
A^\bullet = A^\bullet(V^\bullet) = 
C^\bullet_\Lie(\n^\bullet) \,\,=\,\, (S^\bullet(\n^*[-1]), d)
\]
be the Chevalley-Eilenberg cochain complex of the graded Lie algebra $\n^\bullet$. We consider it
equipped with the total grading, combining the grading induced by that on $\n^\bullet$ and the
cohomological grading by the degree of cochains. With this grading $A^\bullet$ a graded commutative
dg-algebra  situated in degrees $\leq 0$, i.e., it is an object of $\dgAlg_\k$,
and 
\[
A^0\,\,=\,\,  \bigotimes_i S^\bullet (V^i \otimes (V^{i+1})^*)\,\,=\,\,\k \biggl[ \prod_i {\Hom(V^i, V^{i+1})}\biggr]
\]
is the coordinate algebra of the ambient affine space for $\Com(V^\bullet)$. 

\begin{Defi} The derived variety of complexes is the dg-scheme 
\[
R\Com(V^\bullet) \,\,=\,\,\Spec (A^\bullet(V^\bullet)). 
\]
\end{Defi}

The  name ``derived variety of complexes" for $R\Com(V^\bullet)$ is justified by the following fact. 

\begin{prop} (a) The variety of complexes
is identified with the scheme of Maurer-Cartan elements of $\n^\bullet(V^\bullet)$:
\[
\Com(V^\bullet) \,\,=\,\,\bigl\{ D\in\n^1  \, \bigl| \, [D,D]=0 \bigr\}.
\]
(b) We have an identification of commutative algebras 
\[
\k[\Com(V^\bullet)] \,\,\simeq\,\, H^0(A^\bullet) \,\,=\,\, H^0_\Lie(\n^\bullet). \qed
\]
\end{prop}

Additionally, the  name ``derived variety of complexes"
is justified by the following description of the
functor represented by $A^\bullet(V^\bullet)$ on the category
$\dgAlg_\k$. 

\begin{prop}\label{prop:twis-com}
For any $(R,d_R)\in \dgAlg_\k$  the set $\Hom_{ \dgAlg_\k}(A^\bullet(V^\bullet), R)$
is naturally identified with the set of differentials $\delta$ in $R\otimes V^\bullet$  of total degree $+1$ satisfying the
following properties:
\begin{enumerate}
\item[(1)] The Leibniz rule
\[
\delta(rm) = d_R(r) m + (-1)^{\deg(r)} r\delta(m), \quad r\in R, m\in R\otimes V^\bullet.
\]
\item[(2)] The filtration  condition:
\[
\delta(R^i\otimes V^j)  \,\, \subset \,\, \bigoplus_{i'\leq i} R^{i'} \otimes V^{i+j+1-i'}. 
\]

\item[(c)]  $(d_R\otimes 1 +\delta)^2=0$. \qed

\end{enumerate}

\end{prop}

\begin{rems} 
(a) Doubly graded complexes with  the above properties are known as {\em twisted complexes}, cf.  \cite{BK}\cite{FSbook}. 
The construction of $R\Com(V^\bullet)$ can now be said to consist in passing from usual complexes to twisted complexes.

\vskip .2cm

(b) The bidegree (0,1) part of   a differential $\delta$ above, gives rise  to $R$-linear maps 
$\widetilde D_i:  R\otimes V^i\to R\otimes V^{i+1}$. These maps commute with $d$ and so induce
$H^\bullet_d(R)$-linear maps
$D_i: H^\bullet_d(R)\otimes V^i\to H^\bullet_d(R)\otimes V^{i+1}$. Further, the other components of $\delta$  
yield the following:
\begin{itemize}

\item Each  composition $D_{i+1} D_i=0$, and therefore one can speak about the triple  (matrix)  Massey products 
$\langle D_{i+2}, D_{i+1}, D_i\rangle: H^\bullet_d(R)\otimes V^i \to H^\bullet_d(R)\otimes V^{i+2}$.

\item Each $\langle D_{i+2}, D_{i+1}, D_i\rangle$ is equal to 0 modulo indeterminacy, and so one can speak about the
quadruple Massey products  $\langle D_{i+3},  D_{i+2}, D_{i+1}, D_i\rangle: H^\bullet_d(R)\otimes V^i \to H^\bullet_R\otimes V^{i+3}$,
 and so on, so that all the consecutive Massey products
 $\langle D_{i+p}, \cdots,   D_{i+2}, D_{i+1}, D_i\rangle$ are defined and equal to 0 modulo indeterminacy, see \cite{may}.  

\end{itemize}
The cohomology of the dg-algebra $A^\bullet(V^\bullet)$ can be seen therefore  as consisting of ``residual  scalar Massey products":
its elements represent universal  elements of $H^\bullet_d(R)$  which can be constructed whenever we have a datum
exhibiting the existence and vanishing of all the matrix Massey products. 
\end{rems}

We recall, see, e.g., \cite{CFK1} \cite{TVe} that for any derived scheme $X$ over $\k$ and any $\k$-point $x\in X$ we have the tangent
dg-space $T^\bullet_xX$ which is a $\ZZ_{\geq 0}$-graded complex of $\k$-vector spaces defined canonically up up  quasi-isomorphism.

We denote by $\Vect$ the abelian category of $\k$-vector spaces, by $\Cc(\Vect)$ the abelian category of cochain complexes over
 $\Vect$  and by $\Ho (\Vect)$ the homotopy
category of $\Cc(\Vect)$, which is the same as the derived category of  $\Vect$. 

The following property of $R\Com(V^\bullet)$ is analogous to similar properties of other derived moduli schemes \cite{Ka} \cite{CFK1}.

\begin{prop}
Let $D$ be a $\k$-point of $\Com(V^\bullet)$, so that $C^\bullet = (V^\bullet, D)$ is a complex. Then
\[
H^i T^\bullet_D R\Com(V^\bullet) = \begin{cases}
\Hom_{\Cc(\Vect)} (C^\bullet, C^\bullet[1]), & \text{ if }\,\,  i=0, 
\\
\Hom_{\Ho(\Vect)}(C^\bullet, C^\bullet [i+1]), & \text { if }\,\,  i>0. 
\end{cases}
\]
\end{prop}

In the case $i=0$ this reduces to the description of the tangent spaces to the scheme $\Com(V^\bullet)$ given in
\cite[Prop. 17]{cukierman}.

\vskip .2cm

\noindent {\em Proof:}
By definition, $T^\bullet_D R\Com(V^\bullet) = \on{Der}_{A^\bullet} (A^\bullet, \k)$ is the dg-space of derivations $A^\bullet \to \k$ where $\k$ is considered
as a dg-algebra over $A^\bullet$ via evaluation at $D$. Looking at the space $\bigoplus_{i <j}  \Hom_k(V^i, V^j)^*$ of
free generators defining $A^\bullet$, we identify, first as a vector space and then, more precisely,  as a complex: 
\[
\on{Der}_{A^\bullet} (A^\bullet, \k) \,\,= \,\, \bigoplus_{i<j} \Hom_k(V^i, V^j) \,\,=\,\,
\Hom^{\geq 0} (C^\bullet, C^\bullet [1]). 
\]
Here the RHS is the degree $\geq 0$ part of the Hom-complex between $C^\bullet$ and $C^\bullet[1]$.
The cohomology of this part of the Hom-complex is precisely the RHS of the identification claimed in the proposition. 
\qed

 \vskip .2cm

The goal of this paper is the study of the cohomology algebra $H^\bullet(A^\bullet)$, i.e., of the Lie algebra
cohomology of $\n^\bullet$. 

\vskip .3cm

\noindent {\bf B. First examples.} 
Given a $\ZZ_{\leq 0}$-graded commutative dg-algebra $R^\bullet$, we define
\[
H^{\ol 0}(R^\bullet) \,\,=\,\, \bigoplus_i H^{2i}(R^\bullet), \quad H^{\ol 1}(R^\bullet)\,\, =\,\,
\bigoplus_i H^{2i+1}(R^\bullet).
\]
Thus $H^{\ol 0}(R^\bullet)$ is a commutative algebra in the ordinary sense, and $H^{\ol 1}(R^\bullet)$ is
an $H^{\ol 0}(R^\bullet)$-module.  We have the following natural

\begin{que}
Suppose $R^\bullet$ is finitely generated.

\begin{enumerate}
\item[(a)] Is  $H^{\ol 0}(R^\bullet)$ is finitely generated as an algebra? Is $H^{\ol 1}(R^\bullet)$
finitely generated as a module?

\item[(b)] More generally,  what is the algebro-geometric nature of the scheme $\Spec H^{\ol 0}(R^\bullet)$
and of the quasi-coherent sheaf on this scheme corresponding to $H^{\ol 1}(R^\bullet)$?   
\end{enumerate}
\end{que}

It turns out that the answer to the first part is no, and algebras
$A^\bullet(V^\bullet)$ provide easy  counterexamples. We have not seen such examples
in literature.

\begin{ex}
Suppose that $V^1=V^2=V^4=\k$ and all the other $V^i=0$.  We denote by $u,v,w$ the generators of the
1-dimensional vector spaces $\Hom(V^2, V^4)$,  $\Hom(V^1, V^2)$ and $\Hom(V^1, V^4)$. The algebra
$\n$ is then an odd version of the Heisenberg Lie algebra, spanned by $u,v,w$:
\[
\deg(u)=2, \,\,\deg(v)=1, \,\, \deg(w)=3, \quad [u,v]=w, \,\,\, [u,v]=[u,w]=0. 
\]
Let $x,y,z$ be the generators of $A^\bullet$ dual to $u,v,w$, so
\[
\begin{gathered}
A^\bullet = \k[y,z]\otimes\Lambda[x], \quad d(z)=xy,\,\, d(x)=d(y)=0, \\
\deg(x) = -1,\,\, \deg(y)=0, \,\, \deg(z)=-2
\end{gathered}
\]
We find easily that the basis of $H^\bullet(A^\bullet)$ is formed by the classes of $y^a$, $z^b x$ with $a,b\geq 0$,
with the product given by
\[
y^a y^b = y^{a+b}, \quad y^a (z^bx) = 0, \quad (z^bx)(z^{b'}x) = 0. 
\]
That is, 
$H^{\ol 0}(A^\bullet)= \k[y]$ is a polynomial algebra (finitely generated) but $H^{\ol 1}(A^\bullet)$ is
the  infinite-dimensional square-zero ideal
spanned by the $z^bx$, so it is not finitely generated as a module.
In particular, the full cohomology algebra $H^\bullet(A^\bullet)$ is not finitely generated.

\end{ex}

\begin{ex}
Let $V^1=V^2=V^3=V^4 = \k$, with all the other $V^i=0$. Thus $A^\bullet$ is generated by
\[
x_{ij}, \,\, 1\leq i<j \leq 4, \,\,\,\deg(x_{ij})=j-i-1,
\]
with the differential
\[
\begin{gathered}
d(x_{i, i+1})=0,\\
d(x_{13})= x_{12}x_{23}, \,\,\, d(x_{24})=x_{23}x_{34}, \\
d(x_{14}) = x_{12} x_{24} - x_{13}x_{34}.
\end{gathered}
\]
Separating the even and odd parts, we have
\[
\begin{gathered}
A^{\ol 0} = \k[x_{12}, x_{23}, x_{34}, x_{14}] \,\,\oplus \,\, \k[x_{12}, x_{23}, x_{34}, x_{14}] x_{13}x_{24},\\
A^{\ol 1} = \k[x_{12}, x_{23}, x_{34}, x_{14}]x_{13} \,\,\oplus \,\, \k[x_{12}, x_{23}, x_{34}, x_{14}] x_{24}.
\end{gathered}
\]
We write  an elements of $A^{\ol 0}$  as $f+g x_{13}x_{24}$ with $f,g\in \k[x_{12}, x_{23}, x_{34}, x_{14}]$. Then
the differential $d\from A^{\ol 0}\to A^{\ol 1}$ can be written, in the functional notation, as

\[
\begin{gathered}
d( f+g x_{13}x_{24}) \,\,=\,\, (x_{12} x_{24} - x_{13} x_{34}) {\partial f\over\partial x_{14}} +
g x_{12} x_{23} x_{24} - g x_{13} x_{23} x_{34} \,\,= \\
=\,\, \biggl( -x_{34}{\partial f \over\partial x_{14}} - x_{23} x_{34} g \biggr) x_{13} 
+ \biggl( x_{12}{\partial f\over\partial x_{14}} + x_{12} x_{23} g\biggr) x_{24}.
\end{gathered}
\]
Similarly, we write an element of $A^{\ol 1}$ as $\phi x_{13} + \psi x_{24}$, with $\phi, \psi
\in \k[x_{12}, x_{23}, x_{34}, x_{14}]$, and the differential $d\from A^{\overline 1} \to A^{\overline 0}$ as
$$
d(\phi x_{13} + \psi x_{24}) \,\,=\,\, 
\bigl( \phi x_{12}x_{23} + \psi x_{23}x_{34}\bigr) 
-\left(x_{12} {\partial \phi\over\partial x_{14}} +
x_{34} {\partial \psi\over\partial x_{14}} \right) x_{13} x_{24}.
$$

The cohomology ring of $A^\bullet$ is concentrated in even degrees and is infinitely generated as an algebra, with
generators given by the classes of cocycles $x_{12}$, $x_{23}$, $x_{34}$ and $e_i = x_{14}^i x_{23} - i x_{14}^{i - 1} x_{13} x_{24}$,
for $i \ge 1$, and relations $x_{12} x_{23} = x_{23} x_{34} = x_{12} e_i = x_{34} e_i = 0$, $e_i e_j = e_{i + j} x_{23}$.
\end{ex}

\subsection {Representation-theoretic framework.}\label{subsec:rep-fram}

Note that the algebraic group $GL(V^\bullet) = \prod_i GL(V_i)$ acts on  the dg-algebra $A^\bullet$ and thus on the dg-scheme
$R\Com(V^\bullet)$. Let us introduce a representation-theoretic notation in order to analyze the action of
$GL(V^\bullet)$ on $H^\bullet(A^\bullet)$. 

First of all, we can and will assume that the only possible non-zero components of $V^\bullet$
are $V^1, \ldots, V^n$ for some $n\geq 1$. We denote $r_i=\dim(V^i)$ and by  $r=\sum r_i$ the total dimension of $V^\bullet$.
We also denote
\[
\rev=\sum_{i \text{ even}} r_i, \quad \rod=\sum_{i \text{ odd}} r_i, \quad \rev + \rod =r.
\]
We recall that an irreducible representation of $GL(V^{i})$ is given by a
highest weight $\alpha^{(i)}=(\alpha^{(i)}_1 \geq \cdots
\alpha^{(i)}_{r_i})$. We denote the representation corresponding to $\alpha^{(i)}$ by
$\Sigma^{\alpha^{(i)}} (V^i)$ (the Schur functor).

An irreducible representation of $GL(V^\bullet)=\prod GL(V^{(i)})$ is then given 
by a sequence ${[\alpha]} = (\alpha^{(1)}, \ldots, 
\alpha^{(n)})$ of dominant weights for the $GL(V^i)$. We denote by 
\[
\Sigma^{[\alpha]} (V^\bullet) \,\,=\,\, \bigotimes_{i=1}^n \Sigma^{\alpha^{(i)}} (V^i)
\]
the irreducible representation corresponding to $[\alpha]$.  
For  any algebraic representation $U$ of $GL(V^\bullet)$ we denote by
$U_{[\alpha]}\subset U$ the $[\alpha]$-isotypic component of $U$. 

\vskip .2cm

Note that we can regard a sequence $[\alpha]$ of weights for $GL(V^i)=GL(r_i)$ as a  single
weight for $GL(r)$. More precisely, for each  integral weight $\beta = (\beta_1, \ldots,\beta_r)$ of $GL(r)$
we will denote  by
\[
[\beta] \,\,=\,\,\bigl( (\beta_1, \ldots, \beta_{r_1}), (\beta_{r_1+1}, \ldots, \beta_{r_1+r_2}), 
\ldots, (\beta_{r_1+\cdots + r_{n-1}+1}, \ldots, \beta_{r_1+\cdots + r_n})\bigl)
\]
the corresponding sequence of weights of the $GL(r_i)$. We say that $\beta$ is {\em block-dominant},
if each component of $[\beta]$ is a dominant weight for $GL(r_i)$.

In addition, for any dominant integral weight $\beta$ for $GL(r)$ we denote
by $\Sigma^\beta(V)$ the corresponding irreducible representation of $GL(r)$, i.e., the result
of applying the Schur functor $\Sigma^\beta$ to $V$, the total space of the graded vector space $V^\bullet$.

We will be interested in the following. 

\begin{que}\label{que:mult}
Which $\Sigma^{[\alpha]} (V^\bullet)$ appear in $H^\bullet(A^\bullet)$ and with which
multiplicity? In other words, what are the multiplicities of the isotypic components $H^i(A^\bullet)_{[\alpha]}$? 
\end{que}

\subsection{Lie superalgebra framework}
\label{subsec:Lie.framework}

Every $\ZZ$-graded Lie algebra gives after contracting the $\ZZ$-grading to a $\ZZ/2$-grading, a Lie superalgebra.
In particular, for each graded vector space $V^\bullet$ as can consider $\n=\n(V^\bullet)$ as a Lie superalgebra.
Considered for all graded vector spaces $V^\bullet$, the Lie superalgebras $\n(V^\bullet)$
are precisely the nilpotent radicals of all even parabolic subalgebras in all Lie superalgebras  of type $\gl$. 
We recall that a parabolic sub algebra   is called {\em even}, if the corresponding Levi subgroup is a purely
even (non super) reductive algebraic group. In the sequel all the parabolic subalgebras will be assumed even.

Therefore Question \ref{que:mult} can be regarded as a super-analog of the classical question studied by
B. Kostant: determine the cohomology of the nilpotent radical of a parabolic subalgebra
in a semi-simple (or reductive) Lie algebra as a module over the Levi.

\vskip .2cm

However the setup of
parabolic subalgebras is  notationally quite different from that of graded vector spaces, because it is customary 
to consider different parabolic subalgebras
in the same Lie superalgebra $\g$, something that is not natural from the point of view of \S \ref{sec:def}. 
Let us describe this second setup, 
referring to  \cite{M} for general background on Lie supergroups and superalgebras. 

\vskip .2cm

Fix $\mev,  \modd \geq 0$ and  let $m=\mev + \modd$. Consider the Lie superalgebra  $\g = \gl(\mev |\modd)$. A standard realization of $\g$
is as the space  of  $m$ by $m$ matrices $B$  viewed as  2 by 2 block
matrices
\[
B = \begin{pmatrix}
B_{\ev\ev}&B_{\ev\od} 
\\
B_{\od\ev}& B_{\od\od}
\end{pmatrix}   \,\,\in\,\, \on{Mat}(m\times m, \k), \quad B_{\bar i \bar j} \in \on{Mat}(m_{\bar i} \times m_{\bar j}, \k).
\]

Let $\b\subset \g$ be the standard
Borel subalgebra formed by upper triangular matrices $B$. {\em Standard parabolic subalgebras} $\p \supset \b$ are in bijection with
decompositions 
\[
\mev = m_1 + \cdots + m_M, \quad \modd =  m_{M+1} + \cdots +  m_{M+N}, \quad m_i >0. 
\]
Such a decomposition realizes each $B\in \g$ as a block matrix $\|B_{ij}\|$, $i,j=1, \ldots, M+N$, and
\[
\p= \p \bigl(m_1, \ldots, m_M|   m_{M+1}, \ldots m_{M+N} \bigr)  \,\,=\,\,\{ B\mid B_{ij}=0 \text{ unless } i\leq j\}. 
\]
In addition, let $\pi\in S_{M+N}$ be any permutation. We associate to it the  $\pi$-permuted parabolic subalgebra
\[
\p_\pi =  \p_\pi \bigl( m_1, \ldots, m_M|   m_{M+1}, \ldots m_{M+N})  \,\,=\,\,\{ B\mid B_{ij}=0 \text{ unless } \pi(i)\leq \pi(j)\}. 
\]
We denote by
\[
\n_\pi =  \n_\pi \bigl( m_1, \ldots, m_M|   m_{M+1}, \ldots m_{M+N})
\]
the nilpotent radical of $\p_\pi$. 
We will be interested in (even)  parabolic subalgebras up to conjugation (in fact, up to isomorphism). Therefore it will be enough for us
to restrict to the case when $\pi$ is
a $(M,N)$-shuffle as the case of a general $\pi$ is reduced to this one by conjugation, cf. \cite{M}, \S 3.3. 

\begin{exa}
The  class of Borel subalgebras corresponds to the case when all $m_i=1$, so $M=\mev, N=\modd$. For an $(M, N)$-shuffle $\pi$
we denote the corresponding Borel subalgebra simply by $\b_\pi$ and its nilpotent radical by $\n_\pi$. 

\end{exa}

\vskip .2cm

Consider now the setting of \S \ref{sec:def}C, that is, we have a graded vector space $V^\bullet$.
Let us assume that  each $V^i =\k^{r_i}$ is a standard coordinate space and put
\[
V_{\ev}=\bigoplus V^{2i} = \k^{\rev}, \quad V_{\od}=\bigoplus V^{2i+1} = \k^{\rod}.
\]
Then we have an embedding of vector spaces
\begin{equation}\label{eq:n-emb}
\n = \n(V^\bullet) \,\,=\,\, \bigoplus_{i<j} \Hom(V^i, V^j) \,\,\subset \,\, \gl(\rev| \rod). 
\end{equation}
Note that we can consider $\n$ as a Lie superalgera (again, by reducing the $\ZZ$-grading modulo 2). 
Then the embedding  \eqref{eq:n-emb} is a morphism of Lie superalgebras.

Let $m_1, \ldots, m_M$ be all the nonzero $r_{2i}$ written in the order of increasing $i$.
That is,
\[
m_i = r_{\nu_i}, \quad \nu_i \text{ even, } i=1, \ldots, M.
\]
Similarly, let $m_{M+1}, \ldots m_{M+N}$ be all the nonzero $r_{2i+1}$ in the order of increasing $i$. That is,
\[
m_{M+j}  = r_{\nu_{M+j}}, \quad \nu_{M+j} \text{ odd, } j= 1, \ldots, N.
\]
In particular, we have exactly
$M+N$ nonzero spaces $V^i$.  
We have then an $(M,N)$-shuffle permutation $\pi\in S_{M+N}$ which arranges the set  $\{\nu_1, \ldots, \nu_{M+N}\}$
in the increasing order. With these notations, we have

\begin{prop}
The embedding \eqref{eq:n-emb} identifies $\n(V^\bullet)$ with the nilpotent radical  
\[
\n_\pi \bigl(m_1, \ldots, m_M| m_{M+1}, \ldots, m_{M+N}\bigr). \qed
\] 
\end{prop}

Let us also denote by 
\[
\L_\pi = \L_\pi \bigl(m_1, \ldots, m_M| m_{M+1}, \ldots, m_{M+N}\bigr)\,\,=\,\, \prod_{i=1}^{M+N} GL(m_i)
\]
the Levi subgroup corresponding to  $\p_\pi$. As usual, $\L_\pi$ acts on $\n_\pi$ by conjugation, and Question 
\ref{que:mult} is equivalent to the following analog of the question studied by Kostant.

\begin{que}
Which irreducible 
representations of $\L_\pi$ appear in the Lie superalgebra cohomology $H^\bullet(\n_\pi, \k)$ and with which multiplicity? 
\end{que}

\begin{exa}[(Kostant's theorem)]\label{ex:kostant} Suppose $\modd=0$, so $m=\mev$, and therefore each $(M,N)$-shuffle is trivial.
This corresponds to the fact that  
each parabolic subalgebra in $\gl(m)$ is conjugate to a standard one $\p(m_0, \ldots, m_M)$. 
We recall  the classical result of Kostant which describes the action of $\L(m_1, \ldots, m_M)$ on the cohomolofy of $\n(m_1, \ldots, m_M))$
with
coefficients in any irreducible representation of $GL(m)$.

\begin{thm}\label{thm:kostant}
Let $\L=\L(m_1, \ldots, m_M)=\prod_{i=1}^M GL(m_i)$ and $\n = \n(m_1, \ldots, m_M)$.
Then for any integral dominant weight $\beta = (\beta_1 \geq \cdots \geq \beta_m)$
we have an identification of  $\L$-modules
\[
H^p(\n, \Sigma^\beta (\k^m) )\,\,\simeq \,\,\bigoplus_{w\in \on{Sh}(m_1, \ldots, m_M)} 
\Sigma^{[{w(\beta + \rho_m)-\rho_m]}}(\k^m). 
\]
Here $\rho_m = (m-1, m-2, \ldots, 0)$, and $\on{Sh}(m_1, m_2, \ldots, m_M)$
is the set of $(m_1, m_2, \ldots, m_M)$-shuffles inside the symmetric group $S_m$. For each such shuffle
the weight $w(\beta + \rho_r)-\rho_r$ is block-dominant, and 
$\Sigma^{[{w(\beta + \rho_m)-\rho_m]}}(\k^m)$ is the corresponding irreducible representation of $\prod_{i=1}^M GL(m_i)$. \qed
\end{thm}

This has the following interpretation in terms of derived varieties of complexes. 
Suppose the degrees of  all the  nonzero components of $V^\bullet$ have the same parity. 
Let $m_1, \ldots, m_M$ be the dimensions of all the nonzero $V^i$ in the order of increasing $i$.
Then the graded Lie algebra $\n(V^\bullet)$ is identified with the ungraded Lie algebra $\n(m_1, \ldots, m_M)$
as a Lie superalgebra.
Therefore we have an isomorphism of commutative $\ZZ/2$-graded differential algebras
\[
A^\bullet \,\, \simeq \,\, C^\bullet_{\operatorname{Lie}}(\n(m_1, \ldots, m_M), \k)\,\, =\,\,  \Lambda^{\bullet } (\mathfrak{n}^*).
\]
In particular, $A^\bullet$ is finite-dimensional. 
Further, Kostant's theorem (case $\beta=0$)  implies that in this case the action of $GL(V^\bullet)$ on $H^\bullet(A^\bullet)$
has simple spectrum. 
\end{exa}

Our main result is a generalization of this fact.

\begin{thm}\label{thm:simple}
\begin{enumerate}
\item[(a)] For any graded vector space $V^\bullet$ the representation of  $GL(V^\bullet)$ on $H^{\bullet}(A^\bullet)$ has
simple spectrum. 

\item[(b)]  Equivalently, for any $\mev, \modd\geq 0$ and for any (even)  parabolic subalgebra $\p\subset \gl(\mev| \modd)$ with
nilpotent radical $\n$ and  Levi subgroup $\L$,
the action of $\L$ on $H^\bullet(\n, \k)$ has simple spectrum. 
\end{enumerate}
\end{thm}

 \vfill\eject

\section{ Euler characteristic analysis }\label{sec:euler}

\subsection{Formulation of the result and  the  nature of the generating function}\label{subsec:chi-form}

As a first step towards proving Theorem \ref{thm:simple},
in this section we prove the following  Euler characteristic analog.

\begin{thm}\label{thm:euler-simple}
\begin{enumerate}
\item[(a)]  For any graded vector space $V^\bullet$ and any sequence $[\alpha]=(\alpha^{(1)}, \ldots, \alpha^{(n)})$
of Young diagrams,  the Euler characteristic
of the graded vector space 
\[
\Hom_{GL(V^\bullet)}(\Sigma^{[\alpha]}(V^\bullet), H^\bullet(A^\bullet))
\]
is equal to $0$ or $\pm 1$. 

\item[(b)]  Equivalently, for any (even) parabolic subalgebra $\p\subset \gl(\mev| \modd)$ and for any irreducible representation $\Sigma^{[\alpha]}$ 
of  of the corresponding Levi subgroup $\L$,
the Euler characteristic of the graded vector space 
\[
\Hom_{\L}(\Sigma^{[\alpha]}, H^\bullet(\n,\k ))
\]
is equal to $0$ or $\pm 1$. 

\end{enumerate}
\end{thm}

We will prove the formulation (b). Let $\p = \p_\pi = \p_\pi(m_1, \ldots, m_M|m_{M+1}, \ldots, m_{M+N})$, where
$\pi$ is an $(M,N)$-shuffle and $\n=\n_\pi$ be its nilpotent radical. 
For each $i=1, \ldots, M+N$ we introduce a vector of variables $s_i= (s_{i1}, \ldots, s_{i, m_i})$ which are the coordinates
in the maximal torus in $GL(m_i)$
and write
\begin{equation}\label{eq:s-identif}
s= (s_1, \ldots, s_{M+N}) \,\,=\,\, \bigl( (s_{11}, \ldots, s_{1, m_1}), \ldots, (s_{M+N,1}, \ldots, s_{M+N, m_{M+N}})\bigr)
\end{equation}
for the corresponding element of the maximal torus of $\L = \prod_{i=1}^{M+N} GL(m_i)$,  which we denote by $T$.  
Also we write $i \mls{\pi} j$ for $\pi(i) < \pi(j)$.
We will study the generating series for the
(Euler characteristic) character of $T$ on $H^\bullet(\n, \k)$:
\[
F_\pi(s) = \sum_{\lambda\in\ZZ^m} \chi(H^\bullet(\n_\pi, \k)_\lambda) s^\lambda
\]
Here $\lambda$ can be seen as a character of the $T$ and for any $T$-module $U$ we denote
by $U_\lambda$ the $\lambda$-isotypic component of $U$.  

For any Young diagram  $\alpha = (\alpha_1\geq \dots\geq\alpha_m)\in\ZZ^m$
we denote by $\sigma_\alpha(z_1, \ldots z_m)$ the Schur polynomial corresponding to $\alpha$.
Because $\n_\pi$ and $H^\bullet(\n_\pi)$ are $\L$-modules, we can write
\[
F(s) = \sum_{[\alpha] = (\alpha^{(1)}, \ldots, \alpha^{(M+N)}) } c_{[\alpha]} 
\sigma_{\alpha^{(1)}}(s_1) \ldots   \sigma_{\alpha^{(M+N)}}(s_{M+N}),\quad c_{[\alpha]} = 
\chi\bigl( \Hom_{\L}(\Sigma^{[\alpha]}, H^\bullet(\n,\k ))\bigr), 
\]
with the sum being over sequences $[\alpha]$ of Young diagrams.  

\begin{prop}\label{prop:f-phi}
The series $F_\pi(s)$ is an expansion of the  rational function
\[
\Phi_\pi(s) = {\prod\limits_{i < j \atop i,j \le M \text{ or } i,j > M}
\prod\limits_{p=1}^{m_{\mathstrut i}} \prod\limits_{q=1}^{m_{\mathstrut j}}
\left(1 - {s_{ip}\over s_{jq}}\right)\over
\prod\limits_{i \mls{\pi} j \atop i \le M < j \text{ or } j \le M < i}
\prod\limits_{p=1}^{m_{\mathstrut i}} \prod\limits_{q=1}^{m_{\mathstrut j}}
\left(1 - {s_{ip}\over s_{jq}}\right) 
}
\]
in the region 
\[
U_\pi \,\,=\,\,\bigl\{ s \mid \,\, |s_{ip}| < |s_{jq}|  \text{ whenever }  \pi(i) < \pi(j)\bigr\}. 
\]
\end{prop}

\noindent { \em Proof:} Since the Euler characteristic of a complex is equal to the Euler characteristic of its
cohomology, we have
\[
F_\pi(s) = \sum_{\lambda\in\ZZ^m} \chi(C^\bullet(\n_\pi, \k)_\lambda) s^\lambda
\]
Now notice that 
\[
C^\bullet(n_\pi, \k) = S^\bullet(\n_\pi^*[-1]) \,\,=\,\,  S^\bullet ((\n_\pi^{\bar 0})^*) \otimes \Lambda^\bullet( (\n_\pi^{\bar 1})^*),
\]
and 
\[
\n_\pi^{\bar 0} = \bigoplus_{i \mls{\pi} j \atop i \le M < j \text{ or } j \le M < i} \Hom(\k^{m_i}, \k^{m_j}), \quad
\n_\pi^{\bar 1} = \bigoplus_{i < j \atop i,j \le M \text{ or } i,j > M}\Hom(\k^{m_i}, \k^{m_j}).
\]
Here $\k^{m_i}$ is the standard representation of $GL(m_i)$ considered a direct factor of $\L$. Notice further that the
character of $\Hom(\k^{m_i}, \k^{m_j})$ is $\sum_{p=1}^{m_i}\sum_{q=1}^{m_j} s_{jq}/s_{ip}$. Our statement now follows from
the classical formulas for the characters of symmetric and exterior powers (summation of the geometric series). 

\qed

Our proof of Theorem \ref{thm:euler-simple} will be based on the following fact which will allow us to pass between different parabolic
subalgebras of the same $\gl(\mev, \modd)$. 
It shows that in order to obtain the generating series $F_\pi(s)$ for all choices of parabolic subalgebra $\p_\pi$ we need
to understand how the expansion of $\Phi_\pi$ changes when we go from one region to another.

\begin{prop}\label{prop:ratio}
\begin{itemize}
\item[(a)] For any two $(M,N)$-shuffles $\pi, \pi'$ the ratio of the rational functions $\Phi_\pi(s)/\Phi_{\pi'}(s)$
is a Laurent monomial in $s$  taken with coefficient $\pm 1$. 

\item[(b)]
More precisely, 
let $\pi$ and $\pi'$ be two $(M, N)$-shuffles that differ only by one transposition of neighboring elements, say
$i$ and $j$, $i \le M < j$, so that $\pi(j)=\pi (i)+1$ and  therefore $i \mls{\pi} j$ and $j \mls{\pi'} i$. Then
\[
\Phi_{\pi'}(s) = (-1)^{m_i m_j}
{\prod\limits_{p=1}^{m_i} s_{ip}^{m_j} \over \prod\limits_{q=1}^{m_j} s_{jq}^{m_i}}
\Phi_\pi(s).
\]
\end{itemize}
\end{prop}

\noindent {\em Proof:} The numerators of $\Phi_\pi$ and $\Phi_{\pi'}$ coincide. Further, all the factors in the denomiators also coincide except for
the factors $\bigl(1-{s_{ip}\over s_{jq}}\bigr)$ in $\Phi_\pi$ being replaced by the factors $\bigl(1-{s_{jq}\over s_{ip}}\bigr)$ in $\Phi_{\pi'}$.
Our statement now follows by applying the identity
\[
\biggl( 1-{a\over b}\biggr) \biggl/  \biggl( 1-{b\over a}\biggr) \,\,=\,\, -{a\over b}. \qed
\]
We now prove Theorem \ref{thm:euler-simple} in several steps in the rest of this section.

\subsection {Standard Borel}

Let $\p=\b= \b_e$ be the standard Borel subalgebra of upper triangular matrixes in $\gl(\mev, \modd)$.
It corresponds to the $\pi=e$ being the identity permutation. 
In this case $M = \rev$ and $N = \rod$.
To simplify notation we will write $s_i$ instead of $s_{i, 1}$. We have then

\begin{equation}\label{eq:phi-borel}
\Phi (s) = \Phi_e(s)  = {\prod\limits_{1\leq i < j\leq M }\left(1 - {s_i\over s_j}\right) \prod\limits_{M+1 \leq i < j\leq M+N }\left(1 - {s_i\over s_j}\right) \over
\prod\limits_{i\leq M < j} \left(1 - {s_i\over s_j}\right)  },
\end{equation}
and $F(s)=F_e(s)$ is the expansion of this rational function
in the region $|s_i| < |s_j|$ { for }  $i<j$. The two products in the numerators are the Vandermonde determinants:
\[
\prod_{1\leq i < j\leq M}\left(1 - {s_i\over s_j}\right) \,\,= \,\, \sum_{w \in S_{M}} (-1)^{|w|} s^{w(\rho^*_M) - \rho^*_M}, \quad 
\rho_M^* = (0,1, \ldots, M-1),
\]
and similarly for the second product. 
Recall  the Cauchy formula:
\be\label{eq:cauchy-char}
\prod\limits_{i=1}^a \prod\limits_{j=1}^b 
\left(1 - u_i / v_j\right)^{-1} = \sum_{\alpha} \sigma_\alpha(u) \sigma_\alpha(v^{-1}), \quad |u_i| < |v_j|, 
\ee
where the 
sum is taken over all partitions $\alpha = (\alpha_1\geq \cdots \geq \alpha_l\geq 0)$, $l =  \min (a,b)$.  
Recall also the Weyl character formula:
\[
\sigma_\alpha(u_1, \ldots, u_a) \,\,=\,\,
\det(s^{\alpha + \rho_a}) / \det(s^{\rho_a}), \quad \rho_a = (a - 1, \ldots, 1, 0), 
\]
where  
\[
\det(s^\beta) \,\,=\sum_{w\in S_a} (-1)^{|w|} s^{w(\beta)}. 
\]
Combining this together we obtain:
\begin{prop}\label{prop:st1}
The expansion of $\Phi(s)$ in the region $\bigl\{ |s_i|< |s_j| \text{ for }  i<j\bigr\}$
has the form
\[
F(s) =
\sum_\alpha \sum_{w\in S_M \atop v\in S_N} (-1)^{|w|+|v|} (-1)^{M(M-1)/2} (s_1, \ldots, s_M)^{w(\alpha+\rho_M)-\rho_M^*} 
\cdot
(s_{M+1}^{-1}, \ldots s_{M+N}^{-1})^{v(\alpha+\rho_N)-\rho_N}, 
\]
where $\alpha$ runs over partitions  $(\alpha_1\geq \cdots \geq \alpha_l\geq 0)$, $l=  \min (M,N)$. \qed
\end{prop}
This is an explicit Laurent series with coefficients $\pm 1$, and Theorem \ref {thm:euler-simple}
for the standard Borel is proved.

\subsection {Arbitrary Borel in $\gl(N|N)$.}\label{subsec:arbbor}

Let us assume that $M=N$ and $\pi$ is an arbitrary $(N,N)$-shuffle.
We prove Theorem  \ref {thm:euler-simple} for $\p=\b_\pi$.

\sparagraph{Explicit formula for the generating series.}
 In fact we prove a more precise statement identifying
all the coefficients of the series $F_\pi(s)$.

Recall that the roots of $\g$ have the form $\alpha_{ij}= e_i-e_j\in\ZZ^{2N}$ for $i,j\in \{1, \ldots, 2N\}$. 
A root $\alpha_{ij}$ is  odd, if $i\leq M < j \text{ or } j\leq M < j$.
Further, $\alpha_{ij}$ is  positive for $\b_\pi$, if $\pi(i)< \pi(j)$ .
We denote
\[
\delta^\pi = \sum_
{
i\leq M < j 
\atop
{\pi(j) < \pi(i)}
} \alpha_{ji}
\]
the sum of odd roots which are positive for $\b_\pi$ but negative for the standard Borel. 

\vskip .2cm

Let $S$ be the set of bijections $\{1, \ldots, N\} \to \{N+1, \ldots, 2N\}$.  First of all, for each $\phi\in S$
we define the sign factor $(-1)^{N(\phi, \pi)}$ as follows. 
Take $2N$ points 
$a_1, \ldots, a_{2N}$ on the unit circle in $\RR^2$  in the clockwise order but in general position otherwise.
Suppose $\phi\in S$ is given.
For each $i=1, \ldots, N$ we
connect $a_{\pi(i)}$ with $a_{\pi\phi(i)}$ by an arc (straight line interval)
in the disk. We define $N(\phi, \pi)$ be the number of pairwise intersection
points of these arcs.

\vskip .2cm

Next, for each $\phi\in S$ we define the $(\phi, \pi)$-{\em sector} $\Sen(\phi, \pi)\subset\ZZ^{2N}$
to consist of linear integral combinations
\begin{equation}\label{eq:sector}
\sum_{i\in \{1, \ldots, N\}
\atop 
\pi(i) < \pi\phi(i)
}
a_i \alpha_{i, \phi(i)} + \sum_{i\in  \{1, \ldots, N\}
\atop 
\pi(i) > \pi\phi(i)
}
b_i \alpha_{i, \phi(i)}, \quad a_i\geq 0, \, b_i < 0. 
\end{equation}
Note the strict inequalities for the $b_i$. 
Geometrically,  $\Sen(\phi, \pi)$ can be viewed as the set of integer points inside a coordinate
orthant of an $N$-dimensional linear subspace
in $\ZZ^{2N}$. 

For $\lambda\in\ZZ^{2N}$ define
\[
I_\phi^\pi(\lambda) = \begin{cases}
(-1)^{N(\phi, \pi)}, \text { if }  \lambda \in  \Sen(\phi, \pi),\\
0, \text{ otherwise.} 
\end{cases}
\]

\begin{thm}\label{thm:b-expl}
Let $F_\pi(s) = \sum_{\lambda\in\ZZ^{2N}} c^\pi_\lambda s^\lambda$ be the generating series for the Euler characteristics of
the $\lambda$-isotypical components of $H^\bullet(\n_\pi, \k)$. Then:
\begin{enumerate}
\item[(a)] We have
\[
c^\pi_\lambda = \sum_{\phi\in S} I_\pi^\phi(\lambda - \delta^{\pi_1} +\delta^\pi), 
\]
where $\pi_1$ is the shuffle sending $1, 2, \ldots, 2N$ into $1, N+1, 2, N+2, \ldots, N, 2N$. 

\item[(b)] The sum in the RHS of (a) is always equal to $0$ or $\pm 1$. 
\end{enumerate}
\end{thm}

The remainder of this section will be devoted to the proof of Theorem \ref {thm:b-expl}.
We proceed by induction on the length of $\pi$ and on $N$. 

\sparagraph{Base of induction: standard Borel.}  
We first consider the case $\pi=e$ and arbitrary
$N$. In this case the statement is deduced from Proposition \ref{prop:st1} as follows. By putting $\phi = wv^{-1}$ and
$\beta = v(\alpha+\rho_N)$, we rewrite the RHS of the formula in the proposition as follows:
\begin{equation}\label{eq:F-e-1}
(-1)^{N(N-1)/2}  (s_1, \ldots, s_N)^{-\rho_N^*} \cdot (s_{n+1}, \ldots, s_{2N})^{\rho_N}\cdot
\sum_{\beta\in (\ZZ_+^N)_{\neq}} \sum_{\phi\in S} \operatorname{sgn}(\phi) 
\prod_{i=1}^N \left( {s_i\over s_{\phi(i)} } \right)^{\beta_i},
\end{equation}
where $(\ZZ_+^N)_{\neq}$ is the set of integer vectors $(\beta_1, \ldots, \beta_N)$ such that $\beta_i\geq 0$ and $\beta_i\neq \beta_j$ for
$i\neq j$.

Observe that extending the summation in \eqref{eq:F-e-1} from $(\ZZ_+^N)_{\neq}$ to $(\ZZ_+^N)$ does not change the result.
Indeed, if $\beta$ is such that $\beta_i = \beta_j$ for some $i\neq j$, then the monomial
$\prod_{i=1}^N \left( {s_i /s_{\phi(i)} } \right)^{\beta_i}$ is unchanged if we compose $\phi$ with any
permutation preserving $\beta$.  Therefore the sum of such monomials with coefficients $\operatorname{sgn}(\phi)$
will vanish. 

Observe further that for $\pi=e$ the number $N(\phi, e)$ is the number of ``orders" of $\phi$, 
i.e., of pairs $(i<j)$ such that $\phi(i) < \phi(j)$ and so $(-1)^{N(\phi, e)} = (-1)^{N(N-1)/2}\cdot  \operatorname{sgn}(\phi)$. 

We now note that $\delta^e=0$ and  the series
\[
\sum_{\lambda\in\ZZ^{2N}}   \sum_{\phi\in S} I^e_\phi(\lambda - \delta^{\pi_1}) s^\lambda
\]
is precisely the series in
\eqref{eq:F-e-1} in which the summation over $\beta$ is extended to $\ZZ_+^N$. 

This establishes the base of induction. 

\sparagraph {Inductive step: the two expansions.} 
We now assume that Theorem \ref{thm:b-expl} is proved for a given $(N,N)$-shuffle $\pi$ and for
all $(N-1, N-1)$-shuffles $\tau$.  Suppose $|\pi'|=|\pi|+1$ and $\pi'$ and $\pi$ 
differ only by one transposition of neighboring elements, say
$i$ and $j$, as in Proposition  \ref{prop:ratio}(b). We show how to deduce Theorem \ref{thm:b-expl} for $\pi'$.

We recall (Proposition \ref{prop:f-phi}) that $F_\pi(s) =  F_\pi^N(s)$ is the expansion of the rational function
\[
\Phi_\pi(s) =\Phi_\pi^{N}(s) \,\,= \,\,
{\prod\limits_{1\leq i < j\leq N }\left(1 - {s_i\over s_j}\right) \prod\limits_{N+1 \leq i < j\leq 2N }\left(1 - {s_i\over s_j}\right) \over
\prod\limits_{{ i\leq N < j \text{ or } j \leq N < i}
\atop
{\pi(i) < \pi (j) }
} \left(1 - {s_i\over s_j}\right)  }
\]
in the region
\[
U_\pi = \bigl\{  |s_i|< |s_j| \text{ for } \pi(i) < \pi(j)\bigr\}. 
\]
By   inductive assumption,  $F_\pi(s)$ is given by the formula in
Theorem \ref{thm:b-expl}(a).

Let $G_{\pi'}^{N}(s)$ be the Laurent expansion of $\Phi_\pi$ in the region $U_{\pi'}$. 
Because of Proposition \ref{prop:ratio}, $G_{\pi'}^{N}$ is a product of  $F_{\pi'}^{N}$  and 
a monomial with coefficient $\pm 1$, so we start by analyzing the coefficients of $G_{\pi'}$. 

The regions $U_\pi$ and $U_{\pi'}$ are separated by one hyperplane $H= \{s_i=s_j\}$ on which our function $\Phi_\pi(s)$
has a first order pole. To compare $F_\pi$ and $G_{\pi'}$ we need therefore to study how the Laurent expansion
changes after crossing $H$, i.e., after replacing the condition $|s_i/s_j|<1$ to $|s_i/s_j|>1$ while keeping  the magnitudes of
all the other ratios intact. For convenience of the
reader we recall an elementary
1-dimensional situation of such change.

\sparagraph{Laurent re-expansion via the $\delta$-function.} 
Let $\CC[[z, z^{-1}]]$ be the vector space of formal Laurent series $\sum_{n=-\infty}^\infty a_n z^n$, $a_n\in\CC$. We will use the series
\begin{equation}
\delta(z) \,\,=\sum_{n=-\infty}^\infty z^n \,\,\,\in\,\,\, \CC[[z, z^{-1}]]
\
\end{equation}
which  is the Laurent expansion of the delta function at $z=1$. 
Note that
$\CC[[z, z^{-1}]]$ contains the two fields  $\CC((z))$ and $\CC((z^{-1}))$. 
For a rational function $\phi(z)\in \CC(z)$ we denote by $\phi_0\in\CC((z))$ and $\phi_\infty \in \CC((z^{-1}))$ its Laurent expansions
near $0$ and $\infty$. Then, as well known,
\begin{equation}\label{eq:delta-1d}
\left({1\over 1-z}\right)_0  - \left({1\over 1-z}\right)_\infty \,\,=\,\,\delta(z). 
\end{equation}

 \sparagraph{Laurent re-expansion of $\Phi^N_\pi$.}\label{spar:reexp-phi}
 Returning to our situation, let us view $\Phi^{N}_\pi$ as a rational function on $\CC^N\times\CC^N$ with coordinates $s_1, \ldots, s_{2N}$
and consider the space $\CC^{N-1}\times\CC^{N-1}$ with coordinates $s_p, p\neq i,j$. Let
$\varpi: \CC^N\times\CC^N\to \CC^{N-1}\times\CC^{N-1}$ be the projection forgetting $s_i$ and $s_j$. 
We can use $\varpi$ to define pullback $\varpi^*$ on rational functions as well as on formal Laurent series in the $s_p$.
Let $\tau\in S_{2N-2}$ be the $(N-1, N-1)$-shuffle which arranges the set
\[
\pi\bigl( \{ 1, \ldots, 2N\} - \{i,j\}\bigr) \,=\, \pi'\bigl( \{ 1, \ldots, 2N\} - \{i,j\}\bigr)
\] 
in the increasing order. 

\begin{lem}\label{lem:FG}
\begin{itemize}
\item[(a)] We have
\[
F_\pi^{N} -  G_{\pi'}^{N} \,\,=\,\, \bigl(\varpi^* F^{ N-1}_\tau\bigr) \cdot \delta(s_i/s_j).  
\]
\item[(b)] 
We have $F_{\pi'} = -(s_i/s_j) G_{\pi'}$. 
\end{itemize}
\end{lem}

\noindent {\em Proof:} Consider the diagram
\[
\xymatrix{
H\ar[r]^{\hskip -.8cm \epsilon} 
\ar[dr]_\eta
& \CC^N\times\CC^N
\ar[d]^{\varpi}
\\
& \CC^{N-1}\times\CC^{N-1},
}
\]
where $\epsilon$ is the embedding of $H$ and $\eta$ is the composite projection. By inspection of
\eqref{eq:phi-borel}  we note that, first of all, 
$\Phi^{N}_\pi (s) = \Psi(s) \left(1-{s_i\over s_j}\right)^{-1}$, where $\Psi$ does not have pole along $H$.
Further, we note that the restriction  $\epsilon^* \Psi=\Psi|_H$ does not depend on $s_i$ or $s_j$, more precisely,
\[
\epsilon^* \Psi = \eta^* \Phi^{N-1}. 
\]
Part (a) of our lemma follows from this and from \eqref{eq:delta-1d}. Part (b) is a particular case of Proposition
\ref{prop:ratio} (b). 
\qed

\vskip .2cm

We now define $\widetilde{F}_\pi(s) = \widetilde {F}_\pi^N(s)=\sum_\lambda c_\lambda^{\pi'} s^\lambda$,
where $c_\lambda^\pi$ is given by the formula of Theorem \ref{thm:b-expl}(a), and similarly for all
other shuffles and other values of $N$. 
By Lemma \ref{lem:FG} it is enough to establish the following.

\begin{prop}\label{prop:identity}
We have an identity of formal series
\[
\widetilde {F}_\pi^N + (s_j/s_i) \widetilde F_{\pi'}^N\,\,=\,\, 
\bigl(\varpi^* \widetilde F^{ N-1}_\tau\bigr) \cdot \delta(s_i/s_j).  
\]
\end{prop}

\noindent {\em Proof:} We refomulate the required identity in terms of coefficients.
The coefficient of the LHS   at a given $s^\lambda$ is equal to 
\begin{equation}\label{eq:I-phi-sum}
\sum_{\phi\in S}\bigl(  I^\phi_\pi(\lambda-\delta^{\pi_1}+\delta^\pi) + I^\phi_{\pi'}(\lambda-\alpha_{ji}-\delta^{\pi_1}+\delta^{\pi'})
\bigr). 
\end{equation}

The coefficient at $s^\lambda$ in the RHS is zero unless $\lambda$ is of the form $m\alpha_{ij} + \varpi^* \mu$
where $m\in\ZZ$ and  $\mu\in\ZZ^{N-1}\times \ZZ^{N-1}$, in which case it is given by
\begin{equation}\label{eq:psi}
\sum_{\psi\in S'} I^\psi_\tau (\mu - \delta^{\tau_1}+\delta^\tau).
\end{equation}
Here $\tau_1$ is the $(N-1, N-1)$-shuffle defined similarly to $\pi_1$.

We now write the LHS of Proposition \ref{prop:identity} as the sum of two summands:
\[
\sum_{\phi\in S}\bigl(  I^\phi_\pi(\lambda-\delta^{\pi_1}+\delta^\pi) + I^\phi_{\pi'}(\lambda-\alpha_{ji}-\delta^{\pi_1}+\delta^{\pi'})
\bigr) s^\lambda \,\,=\,\,\sum_{\phi(i)=j} + \sum_{\phi(i)\neq j}
\]
The proposition will follow from:

\begin{lem}\label{lem:delta}
(a) 
We have $\sum_{\phi(i)=j} =  \bigl(\varpi^* \widetilde F^{ N-1}_\tau\bigr) \cdot \delta(s_i/s_j)$.

(b) We also have  $\sum_{\phi(i)\neq j}=0$. 
\end{lem}

\noindent {\em Proof:} (a) If $\phi(i)=j$, then
exactly one of the two summands in \eqref {eq:I-phi-sum} is nonzero.
Indeed, looking at the definitions \eqref{eq:sector} of the $(\phi, \pi)$ sector $\Sen(\phi, \pi)$  and the $(\phi, \pi')$-sector
$\Sen(\phi, \pi')$, we see that 
$\alpha_{ij}$ can appear only in the first summand of the condition for belonging to $\Sen(\phi, \pi)$
and in the second summand of the condition for belonging to $\Sen(\phi, \pi')$. Therefore it lies in exactly one of the
two sectors. 

Let 
\[
\psi:  \{1, \ldots, N-1\}   \buildrel \operatorname{mon} \over \longrightarrow     \{1, \ldots, N\} - \{i\} 
\buildrel \phi \over \longrightarrow      \{1, \ldots, N\} - \{j\}   \buildrel \operatorname{mon} \over \longrightarrow   \{ 1, \ldots, N-1\} 
\]
where ``mon" stands for the unique monotone bijection. 

\begin{lem}\label{lem:mu}
The element $\mu\in\ZZ^{N-1}\times\ZZ^{N-1}$ will have $I^\psi_\tau(\mu)$ equal to the nonzero summand in 
\eqref{eq:I-phi-sum}. 
\end{lem}

Indeed, by assumptiion $\pi(j) = \pi(i)+1$  therefore the points $a_i$ and $a_j$ on the circle will be adjacent.
Further, because $\phi(i)=j$, the points $a_i$ and $a_j$ will be joined by an arc, denote it $[a_i, a_j]$. 
Because  $a_i$ and $a_j$ are adjacent, we can assume that $[a_i, a_j]$ does not meet any other arcs $[a_k, a_{\phi(k)}]$.
This means that
$N(\phi, \pi) = N(\psi,\tau)$.    This proves Lemma \ref{lem:mu} and therefore part (a) of Lemma \ref{lem:delta}.

\vskip .2cm

We now prove part (b) of Lemma \ref{lem:delta}.  More precisely, we claim that for $\phi(i)\neq j$, the two summands in
\eqref{eq:I-phi-sum} are either both zero or are the negatives of each other. Indeed, if $\phi(i)\neq j$, then the point $a_i$
on the circle is the endpoint of some arc $[a_p, a_i]$, and $a_j$ is the endpoint of some $[a_q, a_j]$. Interchanging
$i$ and $j$ (to pass from $\pi$ to $\pi'$) will result in $[a_p, a_i]$ being redirected to $a_j$ and $[a_q, a_j]$
redirected to $a_i$. This procedure will change the total intersection number of arcs by 1 modulo 2. 
Lemma \ref{lem:delta} is proved, and therefore Theorem \ref{thm:b-expl} is established.

\subsection {Arbitrary Borel in $\gl (M|N)$.}\label{subsec:arb-borel-chi} 

We want to establish the following.

\begin{prop}
Suppose that Theorem  \ref {thm:euler-simple} holds for all Borel subalgebras in all $\gl(M|M)$. Then it holds for
all Borel subalgebras in all $\gl(M|N)$. 
\end{prop}

\noindent {\em Proof:} By Proposition \ref {prop:ratio}, validity of Theorem  \ref {thm:euler-simple} for all Borel subalgebras
in $\gl(M|N)$ is equivalent to the statement that all Laurent expansions of the function $\Phi=\Phi^{M,N}$ from \eqref{eq:phi-borel}
have all nonzero coefficients $\pm 1$. Suppose 
$M<N$ (the case $M>N$ is treated similarly) and we know
this statement for $\Phi^{M+1, N}$. Let us deduce it for $\Phi^{M,N}$. Denote the extra variable in $\Phi^{M+1, N}$ by $s_0$. Then
\[
\Phi^{M+1, N} (s_0, s_1, \ldots, s_{M+N}) \,\,= \,\, {\prod\limits_{i=1}^M \left( 1-{s_0/s_i}\right ) \over \prod\limits_{j=M+1}^{M+N} \left(1-s_0/s_j\right) }
\cdot
\Phi^{M,N}(s_1, \ldots, s_{M+N}). 
\]
Let $F^{M,N}(s)$ be the Laurent expansion of $\Phi^{M,N}$ in some region $U$, and $F^{M+1, N}$ be the expansion of $\Phi^{M+1, N}$
in the region given by the same inequalities as $U$ together with $|s_0| < |s_p|$, $p=1, \ldots, M+N$. Then
\[
F^{M+1, N} = F^{M,N} \cdot \prod_{i=1}^M \prod_{j=M+1}^{M+N}  \biggl(1 + {s_0\over s_j} + {s_0^2\over s_j^2} + \cdots - 
{s_0\over s_i} - {s_0^2\over s_is_j} - {s_0^3\over s_i s_j^2} -\cdots \biggr). 
\]
This means that $F^{M,N}$ is the sum of all the monomials in $F^{M+1, N}$ which are independent on $s_0$.
So if $F^{M+1, N}$ has all the nonzero coefficients $\pm 1$, then so does $F_{M,N}$. \qed

\subsection {Arbitrary parabolic in $\gl(\mev |\modd)$.}\label{subsec:arb-par-chi}

We now consider the general situation of \S \ref{subsec:chi-form}. That is, we are given decompositions
\[
\mev = m_1 +\cdots + m_M, \,\,\, \modd = m_{M+1} + \cdots + m_{M+N}, \quad m_i \in\ZZ_{>0},
\]
an $(M|N)$-shuffle $\pi$ and denote by $\p_\pi=\p_\pi(m_1, \ldots, m_M| m_{M+1}, \ldots, m_{M+N})$
the corresponding parabolic subgroup and by $\n_\pi$ its nilpotent radical.  
We refine $\pi$ to an $(\mev, \modd)$-shuffle $\Pi$ by replacing each entry $\pi(i)$, $i=1, \ldots, M+N$, 
with the set of integers in the  interval $[m_1 + \cdots + m_{i-1}+1, \, m_1 + \cdots + m_{i}]$  taken in the increasing order. 
We denote by $\b_\Pi$ the Borel subalgebra in $\gl(\mev, \modd)$ corresponding to $\Pi$ and by $\n_\Pi$ its nilpotent
radical. 

\vskip .2cm

Let $F_\pi(s)$ and $F_\Pi(s)$ be the generating series of the Euler characteristic of $T$-isotypic components of
$H^\bullet(\n_\pi, \k)$ and $H^\bullet(\n_\Pi, \k)$ respectively. By Proposition \ref{prop:f-phi}, these series can be
viewed as Laurent expansions of the  two rational functions $\Phi_\pi(s)$ and $\Phi_\Pi(s)$, in the  two  regions $U_\pi, U_\Pi$
respectively. Note that after
identification   \eqref{eq:s-identif}, $U_\Pi = \bigl\{ |s_i|< |s_j|\, \text{ for }  \Pi(i) < \Pi(j)\bigr\}$   
becomes a subset of $U_\pi$, so we can understand both $F_\pi$ and $F_\Pi$ as  the Laurent expansions of $\Phi_\pi$ and  $\Phi_\Pi$
in the same region $U_\Pi$. 

\vskip .2cm

Now, comparing the definitions and identifying the two sets of variables as in 
\eqref{eq:s-identif}, we  see the equality of rational functions
\[
\Phi_\Pi(s) \,\,=\,\,  \Phi_\pi(s) \cdot \prod\limits_{i=1}^{M+N} \prod\limits_{1\leq p< q\leq m_i} \left( 1-{s_{ip}\over s_{iq}}\right)  . 
\]
Since the second factor in the RHS is a Laurent polynomial, this implies the equality of Laurent series
\[
F_\Pi(s) \,\,=\,\,  F_\pi(s) \cdot \prod\limits_{i=1}^{M+N} \prod\limits_{1\leq p< q\leq m_i} \left( 1-{s_{ip}\over s_{iq}}\right)  . 
\]
We now write, as in \S  \ref{subsec:chi-form}, 
\[
F_\pi(s) = \sum_{[\alpha] = (\alpha^{(1)}, \ldots, \alpha^{(M+N)}) } c_{[\alpha]} 
\sigma_{\alpha^{(1)}}(s_1) \ldots   \sigma_{\alpha^{(M+N)}}(s_{M+N}),\quad c_{[\alpha]} \in\ZZ. 
\]
Note that for each $i$  and each Young diagram $\alpha^{(i)}$ the product
\[
P_{\alpha^{(i)}}(s_i) \,\,=\,\,   \sigma_{\alpha^{(i)}} (s_i) \cdot  \prod\limits_{1\leq p< q\leq m_i} \left( 1-{s_{ip}\over s_{iq}}\right)
\]
is a Laurent polynomial with all nonzero coefficients equal to $\pm 1$ (Weyl character formula). 
Further, for different sequences $[\alpha] = (\alpha^{(1)}, \ldots, \alpha^{(M+N)})$ the corresponding product polynomials
\be\label{eq:prod-poly}
P_\alpha(s) = P_{\alpha^{(1)}}(s_1)\ldots P_{\alpha^{(M+N)}}(s_{M+N})
\ee
have disjoint sets of nonzero monomials. Therefore, if some $c_{[\alpha]}\notin \{ 0, \pm 1\}$, then some coefficient
of the Laurent series $F_\Pi(s)$ will be different from $0$ or $\pm 1$. But this is impossible in virtue of the results of
\S \ref{subsec:arb-borel-chi}. Theorem \ref {thm:euler-simple} is now completely proved.

To be more specific, let $\rho^*_{m_i} = (0, \ldots 0, 1, \ldots m_i - 1, 0, \ldots 0)$ be the weights with non-zero components
only in positions from $(m_1 + \cdots + m_{i-1}+2)$ to $(m_1 + \cdots + m_{i})$, for $1 \le i \le M + N$. Assume also that
$\mev > \modd$ and complete the set of variables with $s_{m+1}, \ldots s_{2\mev}$ as in the section \ref{subsec:arb-borel-chi}.
Then we have the following

\begin{cor}
The generating function $F_\pi(s)$ for the Euler characteristic has the following shape
$$
F_\pi(s) = \sum_{\textstyle{\alpha \in Z^{M+N} \atop \alpha_i = 0 \text{ if } i > M+N} \atop \scriptstyle\alpha\text{~--- block-dominant}}
c_{[\alpha]}^\pi \sigma_{\alpha^{(1)}}(s_1) \ldots   \sigma_{\alpha^{(M+N)}}(s_{M+N}),
$$
where coefficients
$$
c_{[\alpha]}^\pi = \sum_{\phi \in S} I_\Pi^\phi \left(\alpha - \delta^{\Pi_1} + \delta^{\Pi} + \sum_{i = 1}^{M + N} \rho^*_{m_i}\right).
$$
\end{cor}

\vfill\eject

\section{Recollections}\label{sec:rec}

\subsection{Mixed complexes and formality}

\noindent {\bf A. Topological motivation.} 
Recall the notion of equivariant formality following \cite{EF}. Let 
 $X$ be a CW-complex endowed with an $S^1$-action. We write $H^\bullet_{S^1} (X) = H^\bullet(X \times_K ES^1, \k)$
for the $S^1$-equivariant cohomology of $X$ with coefficients in $\k$. The fibration $X \times_K ES^1 \to BS^1$ gives rise to the Serre spectral sequence
\begin{equation}\label{eq:SSS}
E_2^{pq} = H^p (BS^1, H^q(X)) \Rightarrow H^{p + q}_{S^1} (X).
\end{equation}

\begin{defi}
We say that space $X$ with the $S^1$-action is {\it equivariantly formal} if the   specral sequence 
\eqref{eq:SSS}
 degenerates at $E_2$.
\end{defi}

We recall from \cite{EF}:
\begin{prop}\label{prop:GM-form}
(a) 
$X$ is equivariantly formal if and only if $H^\bullet_{S^1}(X)$ is a free module over  $H^\bullet(BS^1)=\k[u]$.

(b) If $X$ is equivariantly formal, then $H^\bullet(X) \simeq H^\bullet_{S^1}(X) \big / u H^\bullet_{S^1} (X)$. 
\end{prop}

\vskip 1em

\noindent {\bf B. Mixed complexes.} A {\em mixed complex} is, by definition, a graded $\k$-vector space $C^\bullet$
equipped with two anti-commuting  differentials $d, \delta$ of degrees $(+1)$ and $(-1)$ respectively. See \cite{L}. 
Alternatively,  we may think of a mixed complex $C^\bullet$ as a
graded module over the exterior algebra $\Lambda[d, \delta]$ on two generators:  $d$ of degree $+1$
and $\delta$ of degree $-1$. 

Given two mixed complexes $C^\bullet$ and $D^\bullet$, the tensor product $C^\bullet\otimes_\k D^\bullet$
of graded vector spaces has a natural structure of a mixed complex, with $d$ and $\delta$ extended 
by the Leibniz rule. 

\vskip .2cm

A mixed complex can be seen as a linear algebra analog of a space $X$ with an $S^1$-action. More precisely,
let $X$ be a smooth manifold with a smooth action of $S^1$. Denote by $\theta$ the vector field on $X$ which is the infinitesimal
generator of the action. Then the space of invariant differential forms
\[
\bigl( \Omega^\bullet(X)^{S^1}, \,d = d_{\operatorname{DR}}, \, \delta = i_\theta \bigr)
\]
 equipped with the de Rham differential and the contraction with $\theta$, is a mixed complex. 
 
 Given a mixed complex $(C^\bullet,d, \delta)$, we form a spectral sequence which is an algebraic analog of 
 \eqref{eq:SSS}:
 \begin{equation}\label{eq:SSS2}
\bigl( E_2^{pq} = H^q_d(C^\bullet)^p, d_2 =\delta\bigr) \Rightarrow \HC^{p+q}(C^\bullet).
 \end{equation}
 Here $\HC^\bullet(C^\bullet)$, known as the {\em cyclic cohomology} of $C^\bullet$, is the cohomology of the total complex
 of the {\em Connes double complex} $\Bc C^{\bullet\bullet}$ defined by 
 $(\mathcal B C)^{pq} = C^{q - p}$   with
vertical differential given by $d$ and  the horizontal one given by $\delta$. In fact,  
\eqref{eq:SSS2} is nothing but the standard spectral sequence $('E_n)$ of $\Bc C^{\bullet\bullet}$ but re-graded
by inserting zeroes at all odd levels, 
so that $H^\bullet_d(C^\bullet)$ which normally constitutes $E_1$, appears as $E_2$,
see \cite{L}. 
This motivates the following.

\begin{Defi}
A mixed complex $(C^\bullet, d,\delta)$ is called {\em formal}, if the spectral sequence of the double complex $\Bc C^{\bullet\bullet}$
degenerates at $E_1$ or, what is the same, if the spectral sequence \eqref{eq:SSS2} degenerates at $E_2$. 
\end{Defi}

Note that for the trivial 1-dimensional mixed complex $\k$ (in degree 0)  we have
\begin{equation}
\HC^\bullet (\k) \,=\, \k[u], \quad \deg(u)=2. 
\end{equation}
Therefore for each mixed complex $C^\bullet$ we have that $\HC^\bullet(C^\bullet)$ 
 is naturally a $\k[u]$-module. This action of $u$ comes from its action on Connes'
 double complex, which is the identification $u: \Bc C^{pq} \to \Bc C^{p+1, q+1}$. 
 The following is an algebraic analog of Proposition \ref{prop:GM-form}.

\begin{prop}\label{prop:MC-formal}
(a) A mixed complex $(C^\bullet, d,\delta)$ is {formal} if and only if $\HC^\bullet(C^\bullet)$ is a free $\k[u]$-module. 

(b) If a mixed complex  $(C^\bullet, d,\delta)$ is {formal}, then $H^\bullet_d(C^\bullet) = \HC^\bullet(C^\bullet)/u  \HC^\bullet(C^\bullet)$. 

(c) If $C^\bullet$ and $D^\bullet$ are formal mixed complexes, then so is $C^\bullet\otimes_\k D^\bullet$. 
\end{prop}

\noindent {\em Proof:} 
 If $C^\bullet$ is formal, then we have $\HC^\bullet (C^\bullet) = H^\bullet_d(C^\bullet) \otimes \k[u]$ (tensor product of graded
$\k$-vector spaces). This shows the ``only if" part of (a) as well as (b). Let us now prove the ``if'' part of (a).
Suppose
$\HC^\bullet(C^\bullet)$ is $\k[u]$-free.
 
Let $a\in H^\bullet_d(C^\bullet)$ be an element of the lowest degree with the property that some  higher differential $d_r$
(with the lowest $r$),
applied to $a$, is not trivial. Then $d_r(a)=bu^m$ for some $m$, with $b\in H^\bullet_d(C^\bullet)$. 
This implies that $bu^m$ equals $0$ in $\HC^\bullet(C^\bullet)$. Note further that because of our assumption,
$b$ descends to an element in $\HC^\bullet(C^\bullet)$ which we denote $\bar b$, and we have ${\bar b} u^n=0$.
Because $\HC^\bullet(C^\bullet)$ is $\k[u]$-free, it implies that $\bar b=0$. But this means that some higher differential
is nontrivial on  $b$, which is impossible since it has degree smaller than $a$. 

Finally let us show part (c). Notice that for any mixed complexes $C^\bullet$ and $D^\bullet$
  we have an identification at the  level of Connes double complexes:
\[
\Bc (C^\bullet\otimes_\k D^\bullet) \,\,\simeq \,\, \Bc (C^\bullet) \otimes_{\k[u]} \Bc^\bullet(D^\bullet). 
\]
This gives the   spectral sequence
\[
E_2 = \on{Tor}_\bullet^{\k[u]} \bigl(\HC^\bullet(C^\bullet), \HC^\bullet(D^\bullet)
\bigr)
 \,\,\Rightarrow \,\, \HC^\bullet(C^\bullet\otimes_\k D^\bullet).
\]
 If  $C^\bullet$ and $D^\bullet$  are formal, then the higher $\on{Tor}$s vanish and we have  an identification at the level of cyclic cohomology:
\[
\HC (C^\bullet\otimes_\k D^\bullet) \,\,\simeq \,\, \HC (C^\bullet) \otimes_{\k[u]} \HC^\bullet(D^\bullet),
\]
 which implies that the LHS is free and so $C^\bullet\otimes_\k D^\bullet$ is formal, by (a). 

\qed

\begin{rem}\label{rem:single-gr-SSS}
One can also repackage \eqref{eq:SSS2} into a single graded spectral sequence, i.e., a sequence of  single graded complexes
\be\label{eq:single-SSS}
(E_r^\bullet, \lambda_r), \,\, \deg(\lambda_r)=1-2r,\,\,\, \lambda_r^2=0,
 \quad E^\bullet _{r+1} = H^\bullet_{\lambda_r}(E_r^\bullet),  
\ee
built from the same principles as \eqref{eq:SSS2} but without  introducing several copies of the same space. Explicitly,
\[
E_0^\bullet = C^\bullet, \,\,\, \lambda_0 = d, \,\,\, E_1^\bullet = H^\bullet_d(C^\bullet), \,\,\,\lambda_1 \text{ is induced by } \delta \text{ etc.}
\]
Thus, $C^\bullet$ is formal if and only if all the $\lambda_r$, $r\geq 1$, are zero. 
Note that if some $E_r^\bullet$ is nonzero in only finite range of degrees, then $(E_r)$ converges, and so we have a well defined
term
 $E_\infty^\bullet = E_\infty^\bullet (C^\bullet)$.
  \end{rem}

 \vskip .3cm
 
 \noindent {\bf C. Indecomposable mixed complexes.} Note that mixed complexes, i.e., graded $\Lambda[d,\delta]$-modules,
 form an abelian category.  In particular, we can speak about direct sums,   indecomposables etc. 
 It turns out that classification of mixed complexes is a tame problem of linear algebra.

 \begin{prop}\label{prop:indec-mix}
 Let $C^\bullet$ be an indecomposable mixed complex  which is  a union of finite-dimensional mixed subcomplexes.
   Then $C^\bullet$ is isomorphic to one of the following:
     
\begin{itemize}
\item[(1)] The  4-dimensional ``diamond'' module $D$ (the free rank 1 module over $\Lambda[d,\delta]$):
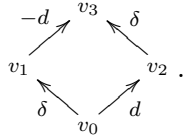
\[
\vcenter{
\def\objectstyle{\scriptstyle}
\def\labelstyle{\scriptstyle}
\xymatrix@-1pc{& v_3 &\\ v_1 \ar^{-d}[ur] && v_2 \ar_{\delta}[ul] \\ & v_0 \ar_{d}[ur] \ar^{\delta}[ul]&}}.
\]

\item[(2)]  Four types of finite-dimensional
 ``zig-zag'' modules $F_n(a, b)$, $n \ge 1$, where $a$ and $b$ are the first and the last arrows respectively in the
diagram
$$
{
 \def\labelstyle{\scriptstyle}
\xymatrix{v_1 \quad \cdots& v_{i-1} \ar^{d}[r] \ar_-{\delta}[l] & v_i & v_{i+1} \ar^-{d}[r] \ar_{\delta}[l] & \cdots \quad v_{N}},
}
$$
and can be either $d$ or $\delta$.  The dimension of the module is $N = 2n$ if $a = b$, and $N = 2n - 1$ if $a \neq b$.

\item [(3)] Five types of infinite-dimensional ``zig-zag'' modules: 2 unbounded on the left side $I_l(d)$ and $I_l(\delta)$
with the last arrow $d$ and $\delta$ respectively, 2 unbounded on the right side $I_r(d)$ and $I_r(\delta)$ with the first
arrow $d$ and $\delta$, and module $I$ unbounded on both sides. We will denote by $v_0$ a generator of the vector space
at the end of the ``zig-zag''.
\end{itemize}

\end{prop}

\noindent {\em Proof:} This follows from a more general result about indecomposable modules over
$\sl(1|1)$ with grading liftable to a $\ZZ$-grading, see \S \ref{subset:gl11}. 

\qed

We can describe formality of a mixed complex in terms of its indecomposable components.

\begin{prop}\label{prop:cohom-indecomp}
\begin{enumerate}[label=\alph*)]

\item[(a)]  The mixed complexes $D$, $I$, $F_n(d, d)$, $I_l(d)$ and $I_r(d)$ are $d$-acyclic.  Thus the spectral sequence
\eqref{eq:SSS2} for any of these modules vanishes, 
   and the modules are formal. 

\item[(b)]  The total $d$-cohomology space
of each of the mixed complexes $F_n(d, \delta)$ and $F_n(\delta, d)$ is one dimensional, generated be the class of $v_{2n}$ and $v_1$
respectively. Hence all the differentials $\lambda_r$, $r\geq 1$,  vanish and the mixed complexes are formal, with 
 \[
 \HC^\bullet (F_n(d, \delta)) \isom \k[u] \cdot [v_{2n}], \quad 
\HC^\bullet(F_n(\delta, d)) \isom \k[u] \cdot [v_{1}].
\]

\item[(c)] The total $d$-cohomology of   each of the the mixed complexes
$I_l(\delta)$ and $I_r(\delta)$ is one dimensional, generated by the class of $v_0$. In particular, these mixed complexes
are formal and $\HC^\bullet (I_*(\delta)) \isom \k[u] \cdot [v_0]$.

\item[(d)]  The total $d$-cohomology of  the mixed complex $F_n(\delta, \delta)$ 
is  two-dimensional, generated by the classes of $v_1$ and $v_{2n}$.
The differential  $\lambda_i = 0$ for $i \neq n$,  while $\lambda_n([v_{2n}]) = [v_1]$. In particular, $F_n(\delta,\delta)$ is not formal,
with
 $\HC^\bullet (F_n(\delta, \delta)) \isom \k[u] / (u^n) \cdot [v_1]$.
\end{enumerate}
\end{prop}

\noindent {\em Proof:} Follows by directly applying the definitions of the higher differentials in the spectral
sequence of a double complex. More precisely, the $r$-th differential in the spectral sequence of $(d,\delta)$
is obtained precisely by considering length $2r$  ``zig-zag" patterns of elements under $d$ and $\delta$. 
\qed

\begin{cor}
Let $C^\bullet$ be  mixed complex  which is the union of finite-dimensional mixed subcomplexes.
Then $C^\bullet$ is formal if and only if  it
does not contain indecomposable summands of type $F_n(\delta, \delta)$. \qed

\end{cor}

\vskip 1em

\noindent{\bf D. Spectral morphisms.} Connes' double complex $\Bc C^{\bullet\bullet}$ of a mixed complex $C^\bullet$ can be
viewed as the polynomial space
$
C^\bullet[u] = C^\bullet\otimes_\k \k[u]
$
with the differentials $d= d\otimes 1$ and $u\delta = \delta\otimes u$. We will also consider the {\em periodized Connes complex}
which is the space of Laurent polynomials
 $
C^\bullet [u, u^{-1}] 
$
equipped with the differential $d+u\delta$.  It is analogous to the construction of  localized equivariant cohomology for
spaces with $S^1$-action and  to the construction of periodic cyclic cohomology of $\k$-algebras \cite{L}. 
Note that the periodized Connes complex is a complex of $\k[u, u^{-1}]$-modules, and
\[
H^\bullet_{d+u\delta} C^\bullet[u, u^{-1}] \,\,=\,\, \HC^\bullet(C^\bullet)[u^{-1}] := \HC^\bullet(C^\bullet) \otimes_{\k[u]}
\k[u, u^{-1}]. 
\]

\begin{Defi}
Let $f: C^\bullet_1\to C^\bullet_2$ be a morphism of mixed complexes. We say that $f$ is a {\em spectral morphism},
if the induced morphism of complexes  of $\k[u, u^{-1}]$-modules
\[
f_*: \bigl(C_1^\bullet[u, u^{-1}], d+u\delta\bigr) \lra \bigl(C_2^\bullet[u, u^{-1}], d+u\delta\bigr)
\]
is a quasi-isomorphism. 
\end{Defi}

\begin{prop}
Let $C_1$ and $C_2$ be formal mixed complexes having
the total $d$-cohomology. 
finite-dimensional.
Let  $f: C^\bullet_1\to C^\bullet_2$ be a spectral morphism. Then
\[
\dim H^{\ev}(C^\bullet_1) = \dim H^\ev(C^\bullet_2), \quad \dim H^{\od}(C^\bullet_1) = \dim H^\od(C^\bullet_2).  
\]
More precisely,  for each of these equalities, there are canonical filtrations in the vector spaces in the left and right hand sides and canonical
isomorphisms between the quotients of these filtrations. 
\end{prop}

Note that the dimensions of individual cohomology spaces $H^i_d(C_1^\bullet)$ and $H^i_d(C_2^\bullet)$
may be different. 

 \noindent {\em Proof:} Using Proposition \ref{prop:MC-formal}(b), we reduce to the following. 
 
 \begin{lem}
 Let $P, Q$ be free finite rank modules over $\k[u]$, and $\Phi\from P\to Q$ be a morphism of modules which becomes
 an isomorphism after  tensoring with $\k[u, u^{-1}]$. Then the fibers $P_0=P/uP$ and $Q_0=Q/uQ$
 have natural filtrations such that  the  coefficients of  the Taylor expansion of $\Phi$ near $u=0$
  define a canonical isomorphism between the associated graded of these two filtrations.
 \end{lem}
 
 \noindent {\em Proof:} We choose trivializations  $P\simeq P_0\otimes k[u]$ and $Q\simeq Q_0\otimes \k[u]$, and then we can speak 
 about the Taylor expansion, i.e., representation of $\Phi$ as a matrix polynomial
  $\Phi(u) = \sum_r \Phi_r u^r$, with $\Phi_r\from P_0\to Q_0$. The
 linear operators $\Phi_r$ do in general, depend on the choice of trivialization. However, $\Phi_0\from P_0\to Q_0$ is canonically defined
 as the map of fibers. Further, $\Phi_1$ descends to
 a canonical map $\phi_1\from \Ker(\Phi_0)\to \Coker (\Phi_0)$. Then, similarly, $\Phi_2$ descends to
   a canonical map $\phi_2\from \Ker(\phi_1)\to\Coker(\phi_1)$ and so on. 
 These maps are canonical because ``each Taylor coefficient is canonically defined modulo the previous ones".
 Alternatively, it is easy to verify the compatibilty of the $\phi_r$ with changes of trivialization.  
 Eventually, for some $r$,  we must have $\Ker(\phi_r)=\Coker(\phi_r)=0$ because $\Phi(u)$ becomes an isomorphism at $u\neq 0$. 
 This means that $\phi_r\from \Ker(\phi_{r-1})\to\Coker(\phi_{r-1})$ is an isomorphism. Further, consider the filtrations
 \[
 \begin{gathered}
 P_0 \supset \Ker(\Phi_0) \supset\Ker(\phi_1) \supset \cdots \supset \Ker(\phi_{r-1}), 
 \\
Q_0 \supset I_{r-1} \supset I_{r-2} \supset \cdots\supset I_{0},  \quad I_\nu := \Ker (Q_0 \to \Coker (\phi_\nu)). 
 \end{gathered} 
 \]
 We see that the maps $\phi_\nu$ define isomorphisms
 \[
 \Ker(\phi_\nu)/\Ker(\phi_{\nu+1}) \lra I_{\nu+1}/I_\nu.   \qed
\]

\begin{rem}
A classical topological situation motivating the concept of a spectral map is the localization theorem in
$S^1$-equivariant cohomology \cite{EF}. Suppose $X$ is an equivariantly formal finite CW-complex with an $S^1$-action,
 and $F=X^{S^1}$ is the fixed point locus. 
 Then the restriction map $H^\bullet_{S^1}(X)\to H^\bullet_{S^1}(F)$ becomes an isomorphism
after tensoring with $\k[u, u^{-1}]$ and and the total dimensions of the usual cohomology spaces of $X$ and $F$ are equal,
while the dimensions of the individual cohomology spaces may be different. 
\end{rem}

 \begin{cor}\label{cor:spectral}
 Let $C_1^\bullet , C_2^\bullet$ be  formal mixed complexes equipped with an action of an algebraic torus $T$
 (commuting with the differentials), and $f\from C_1^\bullet \to C_2^\bullet$ be a spectral morphism compatible with
 the action of $T$.
 Suppose that for each
 character $\chi$ of $T$ the $\chi$-isotypical components of $H^\bullet_d(C_i)$ are finite-dimensional. Then 
 $H^\ev_d(C_1)$ is isomorphic to $H^\ev_d(C_2)$ as a $T$-module, and similarly for $H^\od_d(C_1)$  and  $H^\od_d(C_2)$.
 \qed
 \end{cor}

 %%%%%%%%%%%%%%%%%%%%%%%%%%%%%%%%%%%%%%%%%%%%%%%%%%%%%%%%%%%%%%%%%%

\subsection{Representations and cohomology of $\sl(1 | 1)$} \label{subset:gl11}

The  Lie  superalgebra $\sl(1|1)$ is interesting for us for two reasons. 
  First, mixed complexes
can be seen as particular representations of $\sl(1|1)$. Second, we will use $\sl(1|1)$ to perform
``odd reflections" relating non-conjugate Borel subalgebras in $\gl(\mev, \modd)$. 

\vskip .3cm

\noindent {\bf A. Basics on $\sl(1|1)$.}
We recall, see \cite{UVI} for background,  that the  $\sl(1|1)$ has the basis
\[
H=\begin{pmatrix} 1&0\\0&1
\end{pmatrix},
\quad
Q_+=\begin{pmatrix} 0&0\\1&0
\end{pmatrix},
\quad
Q_-=\begin{pmatrix} 0&1\\0&0
\end{pmatrix}
\]
with $H$ being even, $Q_\pm$ being odd and with the commutation relations
\[
[Q_+, Q_-]=H, \quad [H, Q_\pm]=0. 
\]
In other words, $\sl(1|1)$ is an odd Heisenberg algebra. 
Note that we have  isomorphisms of associative superalgebras
\be\label{eq:mix-sl}
U(\sl(1|1))/ (H) \,\,\simeq\,\, \Lambda[Q_+, Q_-] \,\,\simeq\,\,  \Lambda [d,\delta],
\ee
where the RHS is the  $\ZZ$-graded algebra whose graded modules are mixed complexes. 

\vskip .3cm

\noindent {\bf B. Indecomposables.} 
We now recall the classification of indecomposable representations of $\sl(1|1)$, adapting the results and notation of  \cite{leites},
\cite{germoni}, \cite{UVI}
from the case of $\gl(1|1)$. See also \cite{gotz}. Call a representation of $\sl(1|1)$ {\em liftable}, if its $\ZZ/2$-grading
can be refined to a $\ZZ$-grading so that $\deg(Q_\pm) = \pm 1$. Note that a liftable representation can be extended
to a representation of $\gl(1|1)$, the matrix $\on{diag} (1, -1)$ complementing $\sl(1|1)$ to $\gl(1|1)$,
 acting as the degree operator.

\begin{prop}\label{prop:mod-sl11}
Let $V$ be a liftable indecomposable finite-dimension representation of $\sl(1|1)$. Then, up to possible
change of parity,  $V$ is isomorphic to one
of the following:
\begin{enumerate}

\item[\rm ($\I$)] The trivial 1-dimensional module $\k$. 

\item[\rm ($\II^+$)] The 2-dimensional module $\xymatrix@1{v_1 \ar[r]^{Q_+}&v_2}$ (with $Q_-$ acting by $0$). 

\item[\rm ($\II^-$)] The 2-dimensional module $\xymatrix@1{v_1 &\ar[l]_{Q_-} v_2}$ (with $Q_+$ acting by $0$). 

\item[\rm($\II_\lambda$)] The 2-dimensional module $\xymatrix@1{ v_1 \ar@<.5ex>[r]^{Q_+} & v_2 
\ar@<.5ex>[l] ^{Q_-}
}$, with $Hv_i=\lambda v_i$, $i=1,2$. 

\item[\rm($\III^0$)] The 3-dimensional module 
 $\xymatrix@1{
 v_1&
 v_2 \ar[l]_{Q_-}
  \ar[r]^{Q_+}&v_3
 }$, with $H$ acting by $0$. 

\item[\rm($\III$)]  The diamond module (with $H$ acting by $0$)
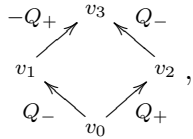
\[
\vcenter{
\def\objectstyle{\scriptstyle}
\def\labelstyle{\scriptstyle}
\xymatrix@-1pc{& v_3 &\\ v_1 \ar^{-Q_+}[ur] && v_2 \ar_{Q_-}[ul] \\ & v_0 \ar_{Q_+}[ur] \ar^{Q_-}[ul]&}}, 
\]

 \item [\rm(Z)] Four types of zig-zag modules
 \[
 \begin{gathered}
 Z^{2p} = \bigl\{ \xymatrix {v_1 & v_2\ar[l]_{Q_-} \ar[r]^{Q_+} &  \,\,\cdots\,\, & \ar[l]_{Q_-} v_{2p-2} \ar[r]^{Q_+} & v_{2p-1}
 & v_{2p}\ar[l]_{Q_-}   }
 \bigr\} , \\
 \ol Z{}^{2p} = \bigl\{
 \xymatrix{
 v_1\ar[r]^{Q_+} & v_2 & \ar[l]_{Q_-}\,\,\cdots  \ar[r]^{Q_+} & v_{2p-2}& \ar[l]_{Q_-} v_{2p-1} \ar[r]^{Q_+}& v_{2p}
 }
 \bigr\},
 \\
 Z{}^{2p+1} = \bigl\{
 \xymatrix{
 v_1\ar[r]^{Q_+} & v_2 & \ar[l]_{Q_-}\,\,\cdots  \ar[r]^{Q_+} & v_{2p-2}& \ar[l]_{Q_-} v_{2p-1} \ar[r]^{Q_+}& v_{2p}&
 \ar[l]_{Q_-} v_{2p+1}
 }
 \bigr\},
 \\
 \ol Z{}^{2p+1} = \bigl\{ \xymatrix {v_1 & v_2\ar[l]_{Q_-} \ar[r]^{Q_+} &  \,\,\cdots\,\, & \ar[l]_{Q_-} v_{2p-2} \ar[r]^{Q_+} & v_{2p-1}
 & v_{2p}\ar[l]_{Q_-} \ar[r]^{Q_+}& v_{2p+1}  }
 \bigr\}. 
  \end{gathered}
 \]
 \end{enumerate}
 Out of these modules, the types {\rm($\I$)}--{\rm($\III$)} are cyclic, while the zig-zag modules are not cyclic except for the following identifications
which exhaust all the isomorphisms of the modules in the list:
 \[
 \I = Z^1 = \ol Z{}^1, \quad \II^+ = \ol Z{}^2, \quad \II^- = Z^2, \quad \III^0 = \ol Z{}^3. 
 \]
\end{prop}

\vskip .3cm

\noindent {\bf C. Multiplication law of liftable cyclic modules.} It follows  \cite{UVI} that the tensor product of two  liftable cyclic
$\sl(1|1)$-modules is again a direct sum of cyclic ones. More precisely, we  have, adapting Prop. 2.3 of {\em loc. cit.} from the case
of $\gl(1|1)$: 

\begin{prop}\label{prop:mult-table}
 Omitting the trivial identifications $I\otimes V=V$, the  multiplication law of liftable cyclic $\sl(1|1)$-modules is the following :
 
\[
\begin{gathered}
 \II_\lambda \otimes \II_{\lambda'} \,\,\simeq \,\, \begin{cases}
 \II_{\lambda+\lambda'} \oplus  \II_{\lambda+\lambda'}, & \text{ if } \lambda+\lambda'\neq 0,\\
 \III, & \text{ if } \lambda+\lambda'=0;
 \end{cases} \\
 \II_\lambda\otimes \II^\pm \,\,\simeq \,\, 2\cdot \II_\lambda;
 \\
 \II^\pm\otimes \II^\pm \,\, \simeq \,\, \II^\pm \oplus \II^\pm, 
 \\
 \II^+ \otimes \II^-  \,\,\simeq\,\, \III; 
 \\
 \II_\lambda\otimes \III^0 \,\, \simeq\,\,  3\cdot \II_\lambda;
 \\
 \II^\pm\otimes \III^0 \,\,\simeq\,\,  \II^\pm \oplus \III;
 \\
 \III^0\otimes \III^0  \,\,\simeq\,\, \II^+ \oplus \III^0 \oplus 2\cdot  \I \oplus \II^-;
 \\
 \II_\lambda\otimes \III \,\, \simeq\,\,  4\cdot \II_\lambda;
 \\
 \II^\pm \otimes \III \,\,\simeq\,\,  2\cdot \III; 
 \\
 \III^0\otimes \III \,\,\simeq\,\,  3\cdot \III;
 \\
 \III \otimes \III \,\,\simeq\,\,  4\cdot \III. 
 \end{gathered} 
\]
\qed
 \end{prop}
 
 \vskip .3cm
 
 \noindent {\bf D. Cohomology with coefficients in cyclic modules.} 
 For any $\sl(1|1)$-module $V$ we denote by $H^i(\sl(1|1), V)$ the $i$-th cohomology of the complex
 of super-vector spaces
 \[
 V\lra \sl(1|1)^*\otimes V \lra \Lambda^2(\sl(1|1)^*) \otimes V \lra \cdots
 \]
 In particular it is iself a super-vector space. 
   We denote by $x,y,t\in \sl(1|1)^*$ the basis dual to
 $Q_+, Q_-, H$. According to the general rules,  $x,y$ become even generators of the Lie algebra cochain complex of $\sl(1|1)$, while $t$
 becomes an odd generator.

 \begin{prop} \label{prop:sl11-coh}
 We have the following identifications:
 \begin{itemize}
 \item[(a)]  $H^\bullet(\sl(1|1), \I) = \k[x,y]/(xy)$ (as an algebra).
 
 \item[(b)] For each  $\lambda\neq 0$ we have $H^\bullet(\sl(1|1), \II_\lambda) = 0$.
 
 \item[(c)] $H^i (\sl(1|1), \III)=
 \begin{cases} \k\cdot [v_3] \text{ (1-dimensional space)}, & \text{ if } i=0,
 \\
 \k\cdot [yv_1+tv_3] \text{ (1-dimensional space)}, & \text{ if } i=1, 
 \\
 0, & \text{ if } i\geq 2.
 \end{cases}
 $
  \end{itemize}
 \end{prop}
 
 \noindent {\em Proof:} (a) By definition, $I=\k$ is the trivial 1-dimensional module, so $C^\bullet(\sl(1|1), \I)= \k[x,y,t]$
 with $x,y$ even, $t$ odd and
 with differential $dt=xy$. So it isomorphic to the complex $\k[x,y] \buildrel xy\over\lra \k[x,y]$
 whence the statement. 
 
 \vskip .2cm
 
 (b) We consider $C^\bullet(\sl(1|1), \II_\lambda)$ as a dg-module over the dg-algebra $C^\bullet(\sl(1|1),\k)=\k[x,y,t]$
 This module is free of rank 2 if we forget the differentials.
 To prove that this module is acyclic, it is enough to  assume $\k$ algebraically closed. Under this assumption,
 we can look at the fibers 
 \[
 F_p\,\,=\,\, C^\bullet(\sl(1|1), \II_\lambda))\otimes_{\k[x,y,t]}  \k_{p}  
 \]
  of $C^\bullet(\sl(1|1), \II_\lambda)$ at various  $\k$-points    $p$ of $\Spec H^\bullet(\sl(1|1), \k)$.
 Each such fiber $F_p$ is a  $\ZZ/2$-graded complex, i.e., a $\ZZ/2$-graded vector space with differential
 $d$ of degree $\ol 1$. The total dimension of $F_p$ is twice the rank of the dg-module, i.e., is equal to 4. 
 
 In order to prove that $C^\bullet(\sl(1|1), \II_\lambda)$ is acyclic, it is enough to prove that $F_p$ is acyclic for
 each $p$. 
   Now, $p$ is a pair $(x_0, y_0)$
 with $x_0, y_0\in\k$ and $x_0 y_0=0$. Let us consider the most degenerate situation $x_0=y_0=0$ (the case
 when $x_0$ or $y_0$ is nonzero being similar).  Assuming $p=(0,0)$, we find that 
 $F_p$
 is identified with the $\ZZ/2$-graded complex
$$
\left\{
\vcenter{
\def\objectstyle{\scriptstyle}
\def\labelstyle{\scriptstyle}
\xymatrix@R=0pt{
\k \cdot v_1 \ar[r]^d & \k\cdot  tv_1\\
\k\cdot t v_2 & \k\cdot v_2 \ar[l]^d
}}
\right\}, \quad d(v_i) = \lambda tv_i, \,\, d(tv_i)=0.
$$
 which is acyclic. 
  
 \vskip .2cm

 To prove (c), let $U$ be the universal enveloping algebra of $\sl(1|1)$, so that for any $\sl(1|1)$-module $M$
 we have $H^\bullet(\sl(1|1),M) = \on{Ext}^\bullet_U(\k, M)$. Note that  $\III$ is identified with $U/(H)$, see \eqref {eq:mix-sl} and so
 has a 2-term resolution
 \be\label{eq:res-III}
 0\to U\buildrel H \over\lra U \to \III\to 0.
 \ee
 
 \begin{lem}
  We have
  \[
  \dim H^i(\sl(1|1), U) = \begin{cases}
  1, & \text{ if } i=1, 
  \\
  0, & \text{ otherwise}. 
  \end{cases}
  \]
 \end{lem}
 
 \noindent {\em Proof of the lemma:} Consider the $U$-module $U[H^{-1}] = U\otimes_{k[H]} \k[H, H^{-1}]$. 
 Considering $\k$ and $U[H^{-1}]$ as $\k[H]$-modules, we see that the groups
 $\Ext^i_{\k[H]}(\k, U[H^{-1}])$ vanish. Since $U$-homomorphisms can be identified with
 $\sl(1|1)$-invariant subspace of $\k[H]$-homomorphisms, we have
 \[
 H^\bullet(\sl(1|1), U[H^{-1}]) \,\,=\,\, \on{Ext}^\bullet_U(\k, U[H^{-1}])=0.
 \]
 On the other hand,  the quotient $U[H^{-1}]/U$ is identified, as an $U$-module,  with the graded
 contragredient module $\Hom_\k (U, \k)$, which implies that
 \[
  \dim H^i(\sl(1|1), U[H^{-1}]/U ) = \begin{cases}
  1, & \text{ if } i=0, 
  \\
  0, & \text{ otherwise}. 
  \end{cases}
  \]
Now, the exact cohomology sequence of the short exact sequence
\[
0\to U \lra U[H^{-1}] \lra U[H^{-1}]/U\to 0
\]
implies the lemma. \qed

We  now prove part (c) of Proposition \ref  {prop:sl11-coh}. 
 The statement about $H^0$, i.e.,  about the space of $\sl(1|1)$-invariants in $\III$, follows at once from
 the definition of the diamond module in Proposition \ref{prop:mod-sl11}. Further, it is immediately to see that
 $yv_1 + tv_3$ is indeed a 1-cocycle of $\sl(1|1)$ with coefficient in $\III$ and that it is not a coboundary. 
 Using the exact cohomology sequence of the resolution \eqref{eq:res-III}, we conclude that there 
 can be no other cohomology. \qed

\vfill\eject

\section{Proof of the Simple Spectrum Theorem}\label{sec:proof-ss}

This section is dedicated to the proof of the main theorem \ref{thm:simple}. The general outline of the proof is similar
to the proof of the Euler characteristic analog \ref{thm:euler-simple}. First we prove the statement for the standard Borel
subalgebra, then via a sequence of reflections extend it to the case of an arbitrary Borel. And finally deduce the result
about arbitrary parabolic subalgebras by relating it to the Borel case by a Hochschild-Serre spectral sequence.

\subsection{Standard Borel}
Let $\g$ be a Lie superalgebra, and $I \subset \g$ be an ideal. Then we have the Hochschild-Serre spectral sequence (cf. \cite{M})
\be\label{eq:HSSS}
E_2^{pq} = H^p(\g / I, H^q(I, \k)) \Rightarrow H^{p + q}(\g, \k).
\ee

We apply it to the case  when $\g = \n = \n_e$ is  the nilpotent radical of the standard Borel subalgebra 
of upper triangular matrices in $\gl(\mev| \modd)$. We denote $V_\ev = \k^{\mev}, V_\od=\k^{\modd}$. 
Further, 
let $I$, $\n_1$, $\n_2$ be the following subalgebras of $\n$:
$$
I  = \left\{ \left(\begin{matrix}0 & B_{\bar 0\bar 1} \\ 0 &0 \end{matrix}\right)\right\}, 
 \quad
\n_1 = \left\{  \left(\begin{matrix}B_{\bar 0 \bar 0} & 0 \\ 0 &0 \end{matrix}\right)\right\}  \cap \n,\quad
\n_2 =\left\{ \left(\begin{matrix}0 & 0 \\ 0 &B_{\bar 1 \bar 1} \end{matrix}\right)\right\}  \cap \n.
$$
Thus  $\n_1$ and $\n_2$ are classical (non-graded) nilpotent Lie algebras, namely the nilpotent radicals of
 the standard Borels in $\gl(\mev)$ and $\gl(\modd)$. Further, $I$ is an odd abelian ideal  in $\n$, identified,
 as an $\n_1\oplus \n_2$-module, with
 $V_\ev \otimes V_\od^*$. 
We may identify the quotient $\n / I$ with the direct sum $\n_1 \oplus \n_2$. Since $I$ is a purely odd abelian Lie superalgebra
we have, applying the Cauchy formula, the following identification of $\n_1\oplus \n_2$-modules:
$$
H^q(I, \k)\,\, \isom\,\, \Lambda^q(I^*) \,\,\isom\,\, S^q(V_\ev^*\otimes V_\od)
\,\, \isom \,\, \bigoplus_{\alpha}\, \Sigma^\alpha V^*_{\bar 0} \otimes \Sigma^\alpha V_{\bar 1}. 
$$
Here the sum is over all partitions $\alpha$, similarly to  \eqref{eq:cauchy-char} where we used the same formula at the level of
characters.

Notice that $\n_1$ acts only on the first factor in the last sum, and $\n_2$ acts only on the second factor, so the cohomology
groups in $E_2$ of the Hochschild-Serre spectral sequence can be written as
$$
H^p(\g / I, H^q(I, \k)) \,\, = \,\, \bigoplus_{\alpha} \bigoplus_{i + j = p}
H^i(\n_1, \Sigma^\alpha V^*_{\bar 0}) \otimes H^j(\n_2, \Sigma^\alpha V_{\bar 1}).
$$
The classical Kostant theorem \ref{thm:kostant} tells us that the weights of $H^i(\n_1, \Sigma^\alpha V^*_{\bar 0})$ 
are precisely  the vertices of the permutahedron $w(\alpha + \rho_\mev) - \rho_\mev$,  and similarly the weights of
$H^j(\n_1, \Sigma^\alpha V_{\bar 1})$ are precisely the  $w(\alpha + \rho_\modd) - \rho_\modd$. Observe that these sets are
disjoint for different $\alpha$. Now the differential $d_r$ in the Hochschild-Serre spectral sequence for $r \ge 2$ goes
from the row $q$ to $q - r + 1$, but $q = |\alpha|$, so the spectral sequence degenerates at $E_2$. And the theorem follows.

\subsection{Arbitrary Borel: outline of the proof}
Suppose $M\geq N$. 
Let $\pi$ be an $(M, N)$-shuffle, in this and the next sections we will establish Theorem \ref{thm:simple} for $\p = \b_\pi$. We proceed
by induction on the length of $\pi$. As before let $\pi'$ be another shuffle that differs from $\pi$ by one transposition
of neighboring elements $i$ and $j$, that is $\pi(j) = \pi(i) + 1$. Assume that theorem holds for $\b_\pi$ and for all
Borel subalgebras in $\gl(M - 1 | N - 1)$. We will show that it also holds for $\b_{\pi'}$.
 Notice that the basis of induction, the case of  $\gl(M - N | 0)$,  follows from the classical
Kostant theorem. 

\sparagraph{Main objects.} 
We start by introducing the main objects that appear in the proof. 
Our situation is that of an ``odd reflection", familiar in representation theory of 
 simple Lie superalgebras, see \cite{Ser1}\cite{Ser2} \cite[\S 3.5]{M}.
 
 \vskip .2cm
 
\begin{figure}[H]
 \centering
 {\hfil\hbox{\begin{tikzpicture}[scale=.5, baseline=(current  bounding  box.center)]
 
 \draw (0,0) -- (11,0);
 \draw (11,0) -- (11,-11); 
 \draw[dashed] (0,0) -- (4, -4);
  \draw[dashed] (11,-11) -- (6, -6);
 \draw (4,-4) --  (11,-4); 
  \draw (4,-6) --  (11,-6); 
    \draw (4,-6) --  (4,0); 
        \draw (6,-6) --  (6,0); 
\node at (5,-5){\small$\sl(1|1)$}; 
\node at (3, -1.5){$J_1$};
\node at (8.5, -1.5){$J_3$}; 
\node at (8.5, -7){$J_2$}; 
\node at (5, -1.5){$L_-$}; 
\node at (8.5, -5){$L_+$}; 
\end{tikzpicture}
}\hfil}
\caption{The structure of $L$.} 
\label{fig:L}
\end{figure}
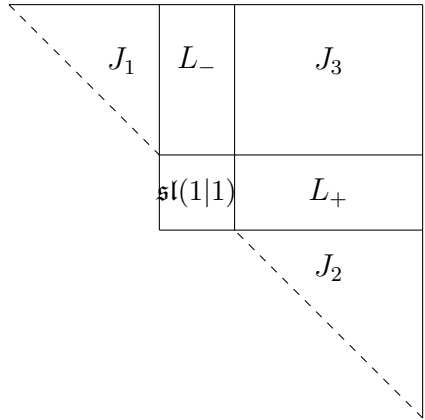

 Consider the Lie superalgebra $L$ generated by the nilpotent radicals $\n_\pi$
and $\n_{\pi'}$ in $\gl(M|N)$. The structure of $L$ can be understood from Figure \ref{fig:L}.
More precisely, this figure represents the conjugate subalgebra $\pi^{-1} L \pi$. 
 
 The direct sum 
 \[
 J_1 \oplus J_2 \oplus J_3 \oplus L_+ \oplus L_- \,\,\,=\,\,\, \n_\pi \cap \n_{\pi'}
 \]
 is an ideal in $L$. The center square  $\sl(1|1)$  is a subalgebra spanned by 
 $e_{ij}, e_{ji}, h = [e_{ij}, e_{ji}]$. The whole Lie superalgebra $L$ is a semidirect sum
  \be\label{eq:semidir}
L = \sl(1|1) \ltimes (\n_\pi \cap \n_{\pi'}). 
\ee
The subalgebra
\[
J \,\,:=\,\,  J_1 \oplus J_2\oplus J_3   \,\,=\,\,  (\n_\pi\cap \n_{\pi'})^{\sl(1|1)}
\]
 is the nilpotent radical of a Borel subalgebra in $\gl(M-1|N-1)$. Note that  $J$ commutes with
$\sl(1|1)$, as indicated by the second equality above.

 We denote the basis of the dual space $\sl(1|1)^*$ by 
\[
t = h^*,\,\, x = e_{ij}^*,\,\,  y = e_{ji}^*,
\] 
so  in the cochain complex $C^\bullet(\sl(1|1))$ we have the relation $dt = xy$. Using
the decomposition \eqref{eq:semidir}, we will view $x,y,t$ as elements of the dual space $L^*$. 

\begin{prop}
 The relation $dt=xy$ holds also in
  $C^\bullet(L)$ as well. 
\end{prop}

\noindent {\em Proof:} This is because
  $h$ is not obtained as a commutator in $L$ in any way essentially different from
$h=[e_{ij}, e_{ji}]$. \qed

\sparagraph{Main steps of the proof.}\label{spar:steps}
 We now present the list of main results that lead to the proof. These results will be proved
one by one later, unless the proof is immediate and presented right away.

\vskip .2cm

Recall that $C^\bullet(L)$ is a commutative dg-algebra, and we denote by $(t,dt)$ the (dg-)ideal in this algebra generated by
$t$ and $dt=xy$. 

\begin{prop}\label{prop:F1}
The natural projection 
\[
q\from C^\bullet(L)\lra C^\bullet(L)/(t,dt)
\]
is a quasi-isomorphism. 
\end{prop}

\begin{prop}\label{prop:F2}
The restriction map $C^\bullet(L)\to C^\bullet(\n_\pi)$ induced by the embedding $\n_\pi\hookrightarrow L$, annihilates the ideal $(t, dt)$
and gives a short exact sequence of complexes
\[
0\to C^\bullet(\n_{\pi'})y \to C^\bullet(L)/(t, dt) \to C^\bullet(\n_\pi)\to 0. 
\]
\end{prop}

\noindent {\em Proof:} This follows immediately by inspecting the Lie superalgebras involved and their cochain complexes. 

\qed
 
 Propositions \ref{prop:F1} and \ref{prop:F2} show that the Lie superalgebra $L$ can be seen as ``interpolating'' between
 $\n_\pi$ and $\n_{\pi'}$, as far as the cohomology is concerned.

\vskip .2cm

Recall that we are studying the {\em spectrum} of various cochain and cohomology spaces, i.e., their decoomposition
under the action of $T=\GG_m^{M+N}$, the maximal torus in $GL(M|N)$. It will be useful to focus in particular on the weight decomposition
with respect to $h\in\on{Lie}(T)$. That is, for any $T$-module $V$ and any $c\in\ZZ$ we denote
\[
V^{h=c} \,\,: = \bigoplus_{\lambda\in T^\vee,\,\, \lambda(h)=c} V_\lambda
\]
the space of $h$-weight $c$. We write $V^h=V^{h=0}$ for the subspace of $h$-weight $0$, i.e., for the kernel of $h$.

\begin{prop}\label{prop:B2}
The embedding $C^\bullet(L)^h\hookrightarrow C^\bullet(L)$ is a quasi-isomorphism. 
\end{prop}

Consider now  the restriction map  from the chain complex
of $L$ to that of $\sl(1|1)\oplus J$, and   denote by
\be\label{eq:res-map}
\res\from C^\bullet(L)^h \lra C^\bullet(\sl(1|1)\oplus J)^h =  C^\bullet(\sl(1|1)\oplus J)
\ee
the induced map on the $h$-weight 0 parts. 

\begin{prop}
The  differential $\delta = i_h$  of degree $-1$ given by the contraction with $h$, makes
(together with the cochain differential $d$),  both $C^\bullet(L)^h$ and $C^\bullet(\sl(1|1)\oplus J)$
into mixed complexes, and $\res$ is a morphism of mixed complexes. 
\end{prop}

\noindent {\it Proof:} In    $C^\bullet(L)$, as in the cochain complex of any Lie algebra, we have the Cartan formula
\[
[d, i_h] \,\,=\,\, \on{Lie}_h,  
\]
where $\on{Lie}_h$ is the operator induced by the coadjoint action of $h$ on $L^*$. So in the $h$-weight 0 part we have
$[d,\delta]=[d, i_h]=0$. \qed

\begin{prop}\label{prop:B3}
Both $C^\bullet(L)^h$ and $C^\bullet(\sl(1|1)\oplus J)$ are formal mixed complexes, and $\rho$ is a spectral morphism.
Therefore (Corollary \ref{cor:spectral}), we have  an isomorphism of $T$-modules
\[
\begin{gathered}
H^\bullet(L) \,=\, H^\bullet(L)^h \,\,\simeq \,\, H^\bullet(\sl(1|1)\oplus J) \,\,=
\\
= \,\, H^\bullet(\sl(1|1))\otimes_\k H^\bullet(J) \,\,=\,\, \k[x,y]/(xy) \otimes_\k H^\bullet(J). 
\end{gathered}
\]
\end{prop}

\begin{rem}
 The previous constructions and results can be seen as ``materialization" of the argument
in \S \ref{subsec:arbbor}\ref{spar:reexp-phi}, more precisely, of Lemma \ref{lem:FG} which compares the Laurent
expansions of the rational function $\Phi^N_\pi$ in two adjacent regions in terms of a function of a smaller number
of variables and a $\delta$-function. More precisely, the cohomology spaces $H^\bullet(\sl(1|1))=\k[x,y]/(xy)$  
``materializes" the
$\delta$-function, since the $T$-weights of this cohomology form an  arithmetic progression infinite on both sides.
The Lie superalgebra $J$, being the nilpotent radical of a Borel in $\gl(M-1|N-1)$, gives, as the generating
series of $H^\bullet(J)$, an expansion of a rational function of smaller number of variables: the analog of the function
$F_\tau^{N-1}$ from Lemma  \ref{lem:FG}. The equality in part (a) of  Lemma  \ref{lem:FG} is now upgraded to the  
short exact sequence of Proposition \ref{prop:F2} combined with the quasi-isomorphism of Proposition \ref{prop:F1}
and the isomorphism of $T$-modules of Proposition \ref{prop:B3}. 
 \end{rem}
 
 We now refine  Proposition \ref{prop:F2} as follows.
  
 \begin{prop}\label{prop:B5}
 The long exact sequence of cohomology associated to the short exact sequence of complexes from Proposition \ref{prop:F2}
 gives:
 \begin{itemize}
 \item[(a)] On the part with $h$-weight $c\neq 0$, isomorphisms of $T$-modules
 \[
 H^\bullet(\n_\pi)^{h=c} \,\simeq\, (H^\bullet(\n_{\pi'})\cdot y)^{h=c}.
 \]
 \item[(b)] On the $h$-weight $0$ part, short exact sequences of complexes
 \[
\xymatrix{
0\to H^{p - 1}(\n_{\pi'})^h \ar^-{\cdot y}[r] & H^p(L) \ar[r] & H^p(\n_\pi)^h\to 0.
}
\]
(Note that $H^\bullet(L)=H^\bullet(L)^h$ by Proposition \ref{prop:B2}.)
 \end{itemize}
 
 \end{prop}
 
 \sparagraph{Completion of the proof modulo the main steps.} Assuming results of n$^\circ$\ref{spar:steps},
 the proof of Theorem  \ref{thm:simple} for $\p = \b_{\pi'}$ is completed as follows.
 
 By inductive assumption, $H^\bullet(\n_\pi)$ has simple spectrum. So the isomorphisms of   Proposition \ref{prop:B5}(a)
 imply that the spectrum of the part  of $H^\bullet(\n_{\pi'})$ with non-zero $h$-weight is also simple. 
 
 Let us concentrate on the part of zero $h$-weight. Recall that $J$ is the nilpotent radical of a Borel   in $\gl(M-1|N-1)$,
 so by our inductive assumtion, $H^\bullet(J)$ has simple spectrum. Denote by $\alpha_x, \alpha_y\in T^\vee$ the
 $T$-weights of $x$ and $y$. Then $\alpha_x=-\alpha_y$ and the subgroup $\ZZ\cdot\alpha_x\subset T^\vee$
 is linearly independent with the subgroup of weights lifted from the maximal torus in $GL(M-1|N-1)$. Therefore
 the isomorphsm of Proposition \ref{prop:B3} implies that $H^\bullet(L)$ has simple spectrum.  Now, the short exact
 sequences from Proposition \ref{prop:B5}(b) imply that $H^\bullet(\n_{\pi'})^h$ has  simple spectrum, since it is a $T$-submodule
 of $H^\bullet(L)$. \qed

 \subsection{Arbitrary Borel: end of proof}
 
 \sparagraph{Proof of Proposition \ref {prop:F1}.}  Let $I=(t,dt)\subset C^\bullet(L)$. Here we prove that
  $q\from C^\bullet(L)\to C^\bullet(L)/I$ is
 a quasi-isomorphism. 
 
 Consider  the filtration of the graded algebra $C^\bullet(L)$ by the powers of the ideal $I$.
The associated graded algebra $\mathrm{gr}_I C^\bullet(L)$ is isomorphic to the tensor product
$C^\bullet(L) / I \otimes S^\bullet(A)$, where $A$ is the acyclic complex $\{t \to dt\}$. All symmetric powers of
$A$ are also acyclic, hence the quotients of the filtration $I^p / I^{p+1}$, for $p \ge 1$ are acyclic as well. This implies that
the map in question is indeed a quasi-isomorphism.

\sparagraph{Proof of Proposition \ref {prop:B2}.}\label{spar:SerreSS}
 Here we prove that the embedding $C^\bullet(L)^h\hookrightarrow C^\bullet(L)$
is a quasi-isomorphism. We start with a general reminder.

Let $\g$ be a Lie superalgebra and $\h\subset\g$ a Lie sub(-super)algebra. Then we have the Hochschild-Serre spectral sequence
\be\label{eq:HSSS-sub}
E_1^{pq} = H^q(\h, \Lambda^p(\g/\h)^*) \,\,\Rightarrow H^{p+q}(\g, \k). 
\ee 
See \cite{fuks} for background. If $\h$ is an ideal, then the $E_{\geq 2}$-part of this spectral sequence is identified with
\eqref{eq:HSSS}. 

We apply this to $\g=L$ and $\h=\sl(1|1)$, so that $\g/\h$ is identified with $\n_\pi \cap \n_{\pi'}$. 
 Let us look at
the structure of $\Lambda^\bullet (\n_\pi \cap \n_{\pi'})^*$ as a representation of $\sl(1|1)$.
Notice that $(\n_\pi \cap \n_{\pi'})^*$ is a direct sum of the $\sl(1|1)$-representations of types
\[
\k = \I, \,\, \k^{1|1} = \II_1, \,\, (\k^{1|1})^* = \II_{-1}.
\]
Therefore all the representations of $\sl(1|1)$ appearing in the sequel (being derived from these three by tensor operations)
 are liftable (possess a $\ZZ$-grading refining the $\ZZ/2$-grading). 
 Using the multiplication table in Proposition
\ref{prop:mult-table} we find that each tensor power and therefore, each exterior power of $(\n_\pi \cap \n_{\pi'})^*$ decomposes into
sum of representations of types $\I$, $\II_\tau$ ($\tau\neq 0$) and $\III$ only.

According to Proposition \ref{prop:sl11-coh}, the cohomology of $\sl(1|1)$ with coefficients in $\II_\tau$, $\tau\neq 0$, 
vanishes. Further, the $h$-weight of the $\sl(1|1)$-cohomology with coefficients in $\I$ and $\III$ is zero.
This follows because $h$ acts trivially on these representations as well as on $\sl(1|1)$ itself (lies in the center).
 Therefore the spectral sequence is trivial in non-zero $h$-weight, which
in turn implies   Proposition \ref {prop:B2}. 

\sparagraph{ Proposition \ref{prop:B3}: formality of $C^\bullet(\sl(1|1)\oplus J)$.}
Here we prove the part of Proposition \ref{prop:B3} which says that $C^\bullet(\sl(1|1)\oplus J)$ is a formal mixed complex.
In fact, we establish a more detailed statement to be used later. We start with a general lemma. 

\begin{lem}\label{lem:HC-as-HKer}
Let $(C^\bullet, d, \delta)$ be a mixed complex, which is free as a $\Lambda[\delta]$-module.
\begin{enumerate}[label=\alph*)]
\item The canonical map $H^p(\Ker \delta) \to \HC^p(C^\bullet)$ is an isomorphism.
\item Let $c \in \Ker \delta$ be a cocycle with respect to $d$. Because of the freeness assumption
there exists $\tilde c \in C^\bullet$, such that $\delta(\tilde c) = c$. Then the cohomology
class of $d\tilde c$ is independent of the choice of $\tilde c$ and under the isomorphism
of (a) is identified with $u \cdot [c] \in \HC^\bullet(C^\bullet)$.
\end{enumerate}
\end{lem}

\noindent{\it Proof:} (a) By the assumption of the lemma the rows of the Connes bicomplex
$\mathcal BC$ are exact, except at the first term. So the corresponding first quadrant
spectral sequence has
$E_1^{pq} = 0$, if $q > 0$ and $E_1^{p0} = \Ker (\delta\from C^p \to C^{p - 1})$, which
implies the statement of the lemma.

(b) Let $\tilde c + \delta c'$ be another lift, then $d(\tilde c + \delta c') - d\tilde c = d \delta c'$ is
a boundary in the complex $\Ker \delta$. The identification with $u \cdot [c]$ follows immediately from
the definition of the Connes bicomplex.

\begin{prop}\label{prop:HCJ}
We have the isomorphism
\[
\HC^\bullet(C^\bullet(\sl(1|1) \oplus J)) \isom \k[x, y] \otimes_\k  H^\bullet(J, \k),
\]
such that the multiplication by $u$ in the cyclic cohomology corresponds to the multiplication by
$xy$ on the right hand side. In particular,  $C^\bullet(\sl(1|1) \oplus J)$ is formal.
\end{prop}

\noindent{\it Proof:} The mixed complex $(C^\bullet(\sl(1|1)\oplus J, d, i_h)$ can be represented as a tensor product of two mixed complexes
(i.e., both $d$ and $\delta$ are extended to the tensor product by the Leibniz rule):
\[
\begin{gathered}
C^\bullet(\sl(1|1)\oplus J) \,\,=\,\, C^\bullet(\sl(1|1))\otimes_\k C^\bullet(J) \,\,= 
\\
=\,\, 
\left(\k[x,y,t], \, d=xy{\partial\over\partial t}, \, \delta = {\partial\over\partial t}\right) \mathop{\otimes}\limits_\k \biggl( C^\bullet(J), \, d, \, \delta =0
\biggr). 
\end{gathered}
\]

So the complex $\Ker \delta$ is isomorphic to $\k[x, y] \otimes C^\bullet(J)$ with $dx = dy = 0$, which implies
the isomorphism claimed in the proposition. Since $C^\bullet(\sl(1|1) \oplus J)$ is a mixed complex with algebra structure it is enough
to show that $u\cdot 1 = xy$, which is clear from $\delta t = 1$ and $dt = xy$. Now, the polynomial ring $\k[x, y]$ is free
as a $\k[xy]$-module, and formality follows from Proposition \ref{prop:MC-formal}(a). \qed

\sparagraph{ $C^\bullet(L)$  as a double complex.}\label{spar:L-double}
 Here we present some preparatory constructions to be used in the proof of the remaining parts of
Proposition  \ref{prop:B3}. 

In the general situation of n$^\circ$ \ref{spar:SerreSS}, assume, in addition, that $\g = \h \ltimes \ken$ is a semi-direct product of
a subalgebra $\h$ and an ideal $\ken$. Then the Hochschild-Serre spectral sequence \eqref{eq:HSSS-sub} comes from the
 fact that  $C^\bullet(\g)$ is realized as the total complex of a double complex. More precisely, 
 \[
 C^\bullet(\g) \,\,=\,\, C^\bullet(\h, C^\bullet(\ken))
 \]
 is nothing but the cochain complex of $\h$ with coefficients in the cochain complex of $\ken$ (on which $\h$ acts). Unravelling the definition,
 we see that 
 \[
 C^\bullet(\gen) = \on{Tot}(C^{\bullet\bullet}, d', d''),  \text{ where } 
 C^{pq} = C^p(\ken) \otimes_\k C^q(\h), \,\,\,  d' =  d_\ken \otimes \on{Id},
\]
and $d_\ken$ is the cochain differential of $\ken$. The differential $d''$ is the cochain differential of $\h$ with coefficients in
$\Lambda^\bullet(\ken^*)$ as an $\h$-module. In particular, if $\h$ is abelian, $d''$ is the standard Koszul differential
corresponding to a family of (super)commuting operators.  

 \vskip .2cm
 
 We apply this to the case $\gen = L$,   with the basis formed by the matrix units $e_{kl}$ and the element $h$.
 We take    $\h = \k\cdot e_{ji}$. It is a 1-dimensional (odd) abelian Lie subalgebra.
  Let also $\ken$ be the subspace in $L$ spanned by all basis elements other than $e_{ji}$.
 Then $\ken$ is an ideal and $\g = \h \ltimes \ken$. 
 We recall that $y$ is the dual basis element corresponding to $e_{ji}$, so $C^\bullet(\h) = \k[y]$ is the 
 usual polynomial algebra. 
 The dg-algebra $C^\bullet(L)$ then contains   $\k[y]$ 
 as a subalgebra and is free as a $\k[y]$-module. The space of free generators is canonically identified with
 $C^\bullet(\ken) = C^\bullet(L)/(y)$.
The above representation of $C^\bullet(\g)$ as a double complex amounts to refining
 the $\ZZ$-grading on $C^\bullet(L)$ to a $\ZZ^2$-grading, by putting
\be\label{eq:double-gr}
C^p(L) = \bigoplus _q C^{p-q, q}, \quad C^{p-q, q} = (C^\bullet(L)/(y))^{p-q} \cdot y^q. 
\ee
Further,  the Koszul differential  corresponding to the action of the 1-dimensional abelian subalgebra $\h$ is 
 the operator induced by the coadjoint action of $e_{ji}$. We denote this operator by $\ell_y = \Lie_{e_{ji}}$. 
So we obtain the following.

\begin{lem}
With respect to the decomposition \eqref{eq:double-gr}, the differential in $C^\bullet(L)$ splits into the sum
$d= \bar d + \ell_y$, where $\bar d$ is the differential in $C^\bullet(L)/(y)$ and $\ell_y = \Lie_{e_{ji}}$.
 In other words, $C^\bullet(L)$ is represented as the total complex of the
double complex $(C^{\bullet\bullet}(L), \bar d, \ell_y)$. \qed
\end{lem}

\sparagraph{ Proposition   \ref{prop:B3}: spectrality of $\res$.} Here we prove that $\res$ is a spectral morphism. 
 For this, we   apply the  double complex realization of n$^\circ$ \ref{spar:L-double} 
  to various subcomplexes in $C^\bullet(L)$.
Furst, we apply it to 
  $h$-weight $0$ part $C^\bullet(L)^h$. Second, we denote by $K$
   the kernel of the restriction morphism $\res$. We note
 that $\res$ is  surjective so that we have a short exact sequence of
mixed complexes
\be
0\to K\lra C^\bullet(L)^h\buildrel\res\over\lra C^\bullet(\sl(1|1) \oplus J) \to 0. 
\ee
The double complex construction applies to $K$ as well. 
In particular, we note the following.

\begin{cor}\label{cor:double}
The complex $K$ is a free module over $\k[y]$ and can be represented 
 as the total complex of the double complex 
$(K^{\bullet\bullet}, \bar d, \ell_y)$, where $K^{p-q,q} = (K^\bullet/(y))^{p-q} \cdot y^q$. \qed
\end{cor}

\begin{lem}\label{lem:HCK-torsion}
The cyclic cohomology $\HC^\bullet(K)$ is a torsion $\k[u]$-module.
\end{lem}

\noindent {\it Proof:} Using the same argument as in Proposition \ref{prop:HCJ} we see that the multiplication by
$u$ in $\HC^\bullet(K)$ coincides with the multiplication by $xy$. Therefore, in order to show that it is a torsion $\k[u]$-module
it is enough to establish that it is a torsion $\k[y]$-module.

Consider the localized complex $K[y^{-1}] = K \otimes_{\k[y]} \k[y, y^{-1}]$. 
By Corollary \ref{cor:double}, we 
 may think of $K[y^{-1}]$ as the total
complex of bicomplex $\tilde K$, with rows $\tilde K^{\bullet q} \isom K / yK$, for all $q \in \ZZ$,
with horizontal differential induced from $K$,
and vertical differential given by the action  of $\ell_y$.

Suppose the columns of this bicomplex are exact with respect to the vertical differential. Then the associated spectral
sequence trivializes at $E_1$, so that cohomology $H^\bullet(K[y^{-1}]) = 0$. This in turn implies vanishing of the
cyclic homology of $K[y^{-1}]$, which can be seen using the spectral sequence \eqref{eq:SSS2}. Hence $\HC^\bullet(K)$
is a torsion $\k[y]$-module.

Exactness of the columns follows from the next lemma.

\begin{lem}
The complex $K / yK$ is free as a $\Lambda[\ell_y]$-module.
\end{lem}

\noindent {\it Proof:} Since $C^\bullet(\sl(1|1) \oplus J)$ is a free $\k[y]$-module we have the short exact sequence:
$$
\xymatrix{K/yK \ar[r] & C^\bullet(L)^h / (y) \ar^-{\res}[r] & C^\bullet(\sl(1|1) \oplus J) / (y)}.
$$
Note that we have an isomorphism of $\Lambda[\ell_y]$-modules
\[
C^\bullet(\sl(1|1) \oplus J) / (y)  \,\,\simeq\,\,  \Lambda^\bullet(\sl(1|1)^*/y) \otimes_\k \Lambda^\bullet(J^*).
\]
Notice further that  $\Lambda^\bullet(J^*)$ has trivial action of $\ell_y$. On the other hand, 
$\sl(1|1)^* / y$ is a free $\Lambda[\ell_y]$-module, and therefore $\Lambda^{>0} (\sl(1|1)^* / y)$ is again a 
free $\Lambda[\ell_y]$-module. This means that $C^\bullet(\sl(1|1) \oplus J) / (y)$ is the direct sum of a free
 $\Lambda[\ell_y]$-module and a trivial module isomorphic to 
 $\Lambda^\bullet(J^*) = \Lambda^0(\sl(1|1)^*/y)\otimes\Lambda^\bullet(J^*)$. 

  We will prove that the $\ell_y$-trivial parts of $C^\bullet(L)^h / (y)$ and of $C^\bullet(\sl(1|1) \oplus J) / (y)$
are identified by the map $\rho$, thus showing that $K/yK=\Ker(\rho)$ is free. 
Indeed,  we have
$$
C^\bullet(L)^h / (y) \ \isom\ \Lambda^\bullet(\sl(1|1)^*/y) \otimes \Lambda^\bullet(L^*_+ \oplus L^*_-)
\otimes \Lambda^\bullet(J^*). 
$$
Note that $L_+ \oplus L_-$  (and so its dual $L_+^* \oplus L_-^*$) 
is a free $\Lambda[\ell_y]$-module. Indeed,  the action on $L_-$ sends $e_{kj}$ to $e_{ki}$, 
identifying the two columns constituting $L_-$ in Fig. \ref {fig:L}.
 Similarly, the action on $L_+$ identifies the two rows constituting $L_+$ in the figure, sending $e_{ik}$ to $e_{jk}$. 
 So the part of $C^\bullet(L)_0 / (y)$ with the
trivial action of $\ell_y$ is isomorphic to 
\[
\Lambda^0(\sl(1|1)^*/y) \otimes \Lambda^0(L^*_+\oplus L^*_-) \otimes\Lambda^\bullet(J^*) \,\,=\,\,
\Lambda^\bullet(J^*)
\]
 and $\res$ maps it isomorphically to the
$\ell_y$-trivial part of $C^\bullet(\sl(1|1) \oplus J) / (y)$. Therefore $K/yK$ is a free $\Lambda[\ell_y]$-module. 

\qed

This establishes Lemma \ref{lem:HCK-torsion} as well. Therefore  $\HC^\bullet(K)[u^{-1}]$, the periodized cyclic cohomology,
vanishes. 
Now the fact that 
 $\res\from C^\bullet(L)^h \to C^\bullet(\sl(1|1) \oplus J)$ is spectral,  follows from the long exact sequence for
cyclic cohomology. 

\sparagraph{Proposition  \ref{prop:B3}: formality of $C^\bullet(L)^h$.} Here we finish the proof of
Proposition  \ref{prop:B3} by showing that $C^\bullet(L)^h$ is a formal mixed complex.

 According to Proposition \ref{prop:MC-formal}(a) we need to show that
$\HC^\bullet(C^\bullet(L)^h)$ is a free $\k[u]$-module, in other words, that multiplication by $u$ has no kernel.
Using the same argument as in the proof of Proposition \ref{prop:HCJ}, we find that multiplication by $u$ is equivalent
to multiplication by $xy$. Therefore it is enough to establish that $x$ and $y$ are not zero divisors in
$\HC^\bullet(C^\bullet(L)^h)$.

First assume that we know that $y$ is not a zero divisor. The the formality follows from the next lemma.

\begin{lem}
Assume that multipliation by $y$ is injective on $\HC^\bullet(C^\bullet(L)^h)$, then multiplication by $x$
is injective as well.
\end{lem}

\noindent{\it Proof:} Consider the following commutative diagram:
$$
\xymatrix{
\HC^\bullet(C^\bullet(L)^h) \ar^-\gamma[r] \ar^{\cdot x}[d]& \HC^\bullet(C^\bullet(L)^h)[y^{-1}] \ar^{\cdot x}[d] \ar^-{\isom}[r]&
\H^\bullet(J) \otimes k[x, y, y^{-1}] \ar^{\cdot x}[d]\\
\HC^\bullet(C^\bullet(L)^h) \ar^-\gamma[r] & \HC^\bullet(C^\bullet(L)^h)[y^{-1}] \ar^-{\isom}[r]& \H^\bullet(J) \otimes k[x, y, y^{-1}]
}
$$
The map $\gamma$ is the canonical map induced by localization, and since $y$ is not a zero divisor, $\gamma$ is injective. 
We have seen before
(see the proof of Lemma \ref{lem:HCK-torsion}) that the cyclic cohomology of the localized complex $K[y^{-1}]$ vanishes, so that
$\HC^\bullet(C^\bullet(L)_0)[y^{-1}]$ is isomorphic to $\HC^\bullet(C^\bullet(\sl(1|1) \oplus J))[y^{-1}]$, which in turn in virtue of
Proposition
\ref{prop:HCJ} is isomophic to $H^\bullet(J) \otimes \k[x, y, y^{-1}]$.

Now, multiplication by $x$ is clearly injective in $H^\bullet(J) \otimes \k[x, y, y^{-1}]$. Combining this with injectivity of $\gamma$
we find that $x$ is not a zero divisor in $\HC^\bullet(C^\bullet(L)^h)$. \qed

Denote by $Z$ the kernel complex $\Ker i_h\from C^\bullet(L)^h \to C^\bullet(L)^h$, which by Lemma \ref{lem:HC-as-HKer} computes the cyclic
cohomology of $C^\bullet(L)^h$. Let us apply the double complex construction of n$^\circ$ \ref{spar:L-double}
  to $Z$ and denote the resulting double complex (refining $Z$)  by $Z^{\bullet\bullet}$ . The
columns of $Z^{\bullet\bullet}$  are isomorphic to $Z / (y)$.

\begin{lem}\label{lem:Z-C}
The canonical projection $C^\bullet(L) \to C^\bullet(\n_\pi)$ restricted to $Z$ induces isomorphism $Z / (y) \isom C^\bullet(\n_\pi)^h$.
\end{lem}

\noindent {\it Proof:} The statement follows immediately from:
$$
Z / (y) = \Ker i_h / (y) \ \isom\ C^\bullet(L)^h / (t, y) \ \isom\ C^\bullet(\n_\pi)^h.\qed
$$

\begin{prop}\label{prop:not-0-div}
The spectral sequence attached to the double complex $Z^{\bullet\bullet}$ degenerates at $E_1$. In particular $y$
is not a zero divisor in $H^\bullet(Z) \isom \HC^\bullet(C^\bullet(L)^h)$.
\end{prop}

\noindent {\it Proof:} Notice that $E_1^{pq} = H^q(C^\bullet(\n_\pi)^h) = H^q(\n_\pi, \k)^h$ for all $p \ge 0$.
Let us fix weight $\lambda$ with respect to the action of the torus $T$ such that $\lambda(h) = 0$,
 and restrict to the $\lambda$-isotypical component. We have
$$
(E^{p\bullet}_1)_\lambda = H^\bullet(\n_\pi, \k)_{\lambda - p \alpha_y} y^p,
$$
where $\alpha_y$ is the root corresponding to $y$.

Now, the explicit formula of Theorem \ref{thm:b-expl} implies that values of the Euler characteristic of
 $H^\bullet(\n_\pi, \k)_{\lambda - p\alpha_y}$,  are either
$\{0, +1\}$ for all $p\geq 0$,  or $\{0, -1\}$ for all $p\geq 0$. By inductive assumption on the $H^\bullet(\n_\pi, \k)$ this implies that the cohomology in weights
$\lambda - p \alpha_y$ is concentrated either in even or odd total degree (here by total degree we understand the degree in
$S^\bullet(\n^*[-1])$). Notice also that the total degree of $y$ is $2$. Therefore the $\lambda$-isotypic part of the spectral sequence can not have non-zero differentials, and since
$\lambda$ was an arbitrary weight satisfying $\lambda(h) = 0$, the proposition follows.

 \qed
 
 This establishes the fact that $t=xy$ is not a zero divisor in $\HC^\bullet(C^\bullet(L)^h)$ and therefore $C^\bullet(L)^h$ is a formal mixed
 complex, thus completing the proof of Proposition   \ref{prop:B3}. 

 \sparagraph{Proof of Proposition \ref{prop:B5}.} We include the short exact sequence of Proposition 
 \ref{prop:F2} into a commutative diagram of short exact sequences
 
\begin{equation}\label{eq:RLES-diag}
\vcenter{
\xymatrix{0\ar[r]&
\mathcal (t,y) \cdot C^\bullet(L) \ar[r] \ar[d]^{q'} & \C^\bullet(L) \ar[r] \ar[d]^{q} & \C^\bullet(L)/(t,y)  \ar[d]^{=}\ar[r]&0\\
0\ar[r]&  y \cdot C^\bullet(\n_{\pi'})  \ar[r] & \C^\bullet(L) / (t, dt) \ar[r] & \C^\bullet(\n_\pi)\ar[r]&0
}}
\end{equation}
Here the vertical maps are the canonical projections, $q$ is the quasi-isomorphism of Proposition \ref{prop:B2}, the right vertical
arrow is the identity, and $q'$, the restriction of $q$,  is  a quasi-isomorphism since the two other vertical arrows are. 

Proposition \ref{prop:B5} is a statement about the long exact sequence (LES) of cohomology associated to the bottom row in
\eqref{eq:RLES-diag}. Part (a) follows at once from Proposition \ref{prop:B2}: since for $c\neq 0$
\[
H^p\bigl(C^\bullet(L)/(t,dt)\bigr)^{h=c} \,\,=\,\, H^p(C^\bullet(L))^{h=c} \,\,=\,\,0, \quad p\in \ZZ,
\]
the cobondary maps in the long exact sequence of cohomology are the required isomorphisms.

Let us prove part (b). For this we look at the $h$-weight $0$ part of the entire diagram \eqref{eq:RLES-diag}. 
We need to prove that the coboundary maps in the LES associated to the bottom row are zero. Because the
vertical arrows are quasi-isomorphisms, this is
equivalent to showing that the coboundary maps in the LES associated to the top row are zero. 
By Lemma \ref {lem:Z-C}, the top row contains a short exact sub-sequence
\be\label{eq:subsequence}
0\to yZ \lra Z \lra Z/yZ\to 0,
\ee
with the same rightmost term $Z/yZ = C^\bullet(\n_\pi)^h$. We view this subsequence as a partial splitting of the
top row of the $h=0$ part of \eqref{eq:RLES-diag} and conclude that
  the coboundary maps associated to that top row
factor through the coboundary maps associated to \eqref{eq:subsequence}. 
But Proposition \ref{prop:not-0-div} means that the maps
$
H^p(yZ) \to H^p(Z)
$
are injective, and so the coboundary maps $H^{p-1}(Z/yZ)\to H^p(yZ)$ are zero. This completes the proof.

\subsection{Arbitrary parabolic}

Our argument will be a refinement of that in \S \ref{subsec:arb-par-chi}, but to visualize the situation, we switch back
 to the (equivalent) notation of \S \ref{subsec:rep-fram}. That is, we denote by $r_1, \ldots, r_n$ the dimensions of the vector
 spaces $V^i = \k^{r_i}$ (some of which may be $0$).  We consider the Lie superalgebra
   $\n = \n(V^\bullet)=\bigoplus_{i<j} \Hom(V^i, V^j)$ formed by block upper triangular 
 $r$ by $r$ matrices, with the parity of the summands induced by the parity of $i$ and $j$.  
 We will prove that the $GL(V^\bullet)$-spectrum of $H^\bullet(\n)$ is simple.
 Note that $\n$ is embedded into the Lie superalgebra $\m$ formed by matrices which are strictly upper triangular in the usual (not block) sense.
 This $\m$ is the nilpotent radical of a Borel subalgebra $\b$ in $\gl(r_{\ev}| r_{\od})$ containing $\n$, and we know that $T$-spectrum of $H^\bullet(\m)$ is simple.
 \vskip .2cm

Let $\u_m$ denote the nilpotent radical of the standard Borel in $\gl(m)$, i.e., the usual Lie algebra of strictly upper triangular matrices.  
In this notation, we have a semidirect sum decomposition
\[
\m = \u \ltimes \n, \quad \u := \u_{r_1}\oplus \cdots \oplus \u_{r_n}
\] 
and the corresponding spectral sequence
\be\label{eq:SS-4}
E_2^{pq} = H^p(\u, H^q(\n)) \Rightarrow H^{p+q}(\b). 
\ee
We identify, as modules over $GL(V^\bullet)$,
\[
H^q(\n) = \bigoplus_{[\alpha]}  \Sigma^{[\alpha]}(V^\bullet) ^{\oplus c_{[\alpha]}^p}, 
\]
with some yet unknown multiplicities $c_{[\alpha]}^p$. Here $[\alpha]= (\alpha^{(1)}, \ldots, \alpha^{(n)})$ runs over sequences of
dominant weights for the $GL(V^i)$. By the classical Kostant theorem, we have an identification of $T$-modules
\be\label{eq:ide-T-mod}
H^p(\u, \Sigma^{[\alpha]}(V^\bullet)) \,\,\simeq\,\, \bigoplus_{{w_1\in S_{r_1}, \ldots, w_n\in S_{r_n}}
\atop {l(w_1)+\cdots + l(w_n)=p} }\bigotimes_{i=1}^n\,\,  \k_{w_i(\alpha^{(i)} + \rho^{(i)})
-\rho^{(i)}}
\ee
Here we identify $T=\prod_{i=1}^n \GG_m^{r_i}$ and for a character $\lambda\in (\GG_m^{r_i})^\vee$ write $\k_\lambda$ for the corresponding
1-dimensional representation. 

Now,  it follows from the observation of  \S \ref{subsec:arb-par-chi}, for different sequences $[\alpha]$ the supports of the $T$-modules
given by the RHS of \eqref{eq:ide-T-mod} do not intersect: their characters are given by the product polynomials $P_\alpha(s)$ from 
 \eqref{eq:prod-poly}. Therefore the only possibility for the differentials of the spectral sequence \eqref{eq:SS-4} is to act between the
 spaces \eqref{eq:ide-T-mod} with the same $\alpha$ and different $p$. 
 Each differential increases $p$ by $2$ or more and is $T$-equivariant. Now, for different $p$,  the subsets
 \[
\Lambda_p = \bigl \{ w= (w_1, \ldots, w_n) \in S_{r_1} \times \cdots\times S_{r_n} \mid  \sum l(w_i) = p\bigr\}
 \]
 are, obviously, disjoint. Therefore applying $w$ from different $\Lambda_p$ to the sequence of strictly
 dominant weights $\alpha^{(i)}+\rho^{(i)}$, we will never  get the same weights.
  Therefore all  the differentials in the spectral sequence are trivial. It further follows that if $c_{[\alpha]}^p \>1$ for some $p$, then
 there will be a $T$-isotypic component of $H^\bullet(\m)$ with multiplicity $>1$ which is impossible. This finishes the proof. 
 
\vfill\eject

\let\thefootnote\relax\footnote {
M.K.: Kavli Institute for Physics and Mathematics of the Universe (WPI), 5-1-5 Kashiwanoha, Kashiwa-shi, Chiba, 277-8583, Japan.
Email: {mikhail.kapranov@ipmu.jp}

\vskip .2cm

 S.P.: \newblock {Department of Mathematics, Yale University, New Haven CT 06520 USA. } \break
 Email:   {svyatoslav.pimenov@yale.edu}

}

\end{document}